\documentclass[12pt,a4paper]{amsart}
\usepackage[T1]{fontenc}
\usepackage[utf8]{inputenc}
\usepackage{lmodern}
\usepackage[margin=2.5cm]{geometry}
\usepackage{amsmath,amssymb,amsthm,mathtools}
\usepackage[only,llbracket,rrbracket]{stmaryrd}%
\usepackage{booktabs}
\usepackage{tabularx}
\usepackage{graphicx}
\usepackage{pgfplots}
\pgfplotsset{compat=1.18}
\usepackage{enumitem}
\usepackage{xurl}
\usepackage{algorithm}
\usepackage{algpseudocode}
\usepackage[hidelinks]{hyperref}
\usepackage[capitalise,noabbrev]{cleveref}
\usepackage{color}%

\newtheorem{theorem}{Theorem}[section]
\newtheorem{lemma}[theorem]{Lemma}
\newtheorem{proposition}[theorem]{Proposition}
\newtheorem{corollary}[theorem]{Corollary}
\theoremstyle{definition}
\newtheorem{definition}[theorem]{Definition}
\newtheorem{example}[theorem]{Example}
\theoremstyle{remark}
\newtheorem{remark}[theorem]{Remark}
\numberwithin{equation}{section}

\DeclareMathOperator{\diam}{diam}
\newcommand{\norm}[1]{\lVert #1\rVert}
\newcommand{\abs}[1]{\lvert #1\rvert}
\newcommand{\jump}[1]{\llbracket #1\rrbracket}
\newcommand{\T}{\mathcal{T}}
\newcommand{\N}{\mathcal{N}}

\newcommand{\OO}{\mathcal{O}}
\newcommand{\R}{\mathbb{R}}
\newcommand{\Id}{I}
\newcommand{\Rec}{\mathcal{R}}
\newcommand{\tol}{\mathrm{tol}}

\newcommand{\PiER}{\Pi_{\mathrm{ER}}}
\newcommand{\PiERf}{\Pi_{\mathrm{ER}}^{f}}
\newcommand{\PiC}{\Pi_{\mathrm{C}}}

\definecolor{brightcerulean}{rgb}{0.11, 0.67, 0.84} %
\definecolor{cardinal}{rgb}{0.77, 0.12, 0.23} %
\definecolor{forestgreen}{rgb}{0.13, 0.55, 0.13}

\definecolor{lightgray}{gray}{0.7}%

\title[Mesh transfers and time error estimation]{On mesh transfers and time error estimation in adaptive finite element methods for parabolic problems}
\author{Stefan Karch}
\address{TU Wien, Institute of Analysis and Scientific Computing, Wiedner Hauptstr. 8-10, 1040 Wien, Austria}
\email{stefan.karch@asc.tuwien.ac.at}

\begin{document}

\begin{abstract}
In space and time adaptive finite element methods, the discrete solution is transferred onto a different mesh whenever the spatial mesh changes.
Refinement is commonly regarded as harmless, since nodal interpolation then reproduces the discrete solution exactly.
However, we show that for the implicit Euler method applied to the heat equation, the transferred solution and time derivative, in general, violate the discrete equation on the new mesh.
Unless the resulting transfer residual vanishes, the corresponding time error estimator overestimates the local error, for small time steps, by a factor inversely proportional to the time step size.
To address this issue, we propose transfers based on the elliptic reconstruction, which significantly reduce the transfer residual or even remove it. 
In numerical experiments, we confirm the analysis and demonstrate the advantages of the proposed transfers in a space and time adaptive example.
\end{abstract}

\maketitle

\section{Introduction}\label{sec:intro}

In this work, we study how the transfer of a discrete solution to a new finite element space affects time error estimation in space and time adaptive finite element algorithms for parabolic equations.
In particular, we focus on mesh refinement.
We show that nodal interpolation, which is commonly used after refinement, preserves the discrete solution but perturbs its discrete time derivative.
This perturbation enters the time error estimator and thus leads to an overestimation of the local time error by a factor proportional to the inverse time step size.
To address this issue, we propose mesh transfers based on the elliptic reconstruction and demonstrate that these avoid such perturbations.

In space and time adaptive computations, the spatial mesh is changed every few time steps, and the discrete solution has to be transferred to the new finite element space in a suitable way.
Coarsening is commonly regarded as the critical part of a mesh change, since the discrete solution must be approximated in a space that cannot represent it exactly, and repeated uncontrolled mesh changes can even prevent convergence~\cite{Dup82}.
Under refinement, however, the old finite element space is contained in the new one, and nodal interpolation represents the discrete solution exactly.
One may therefore expect that a time error estimator based on the discrete time derivatives remains unaffected by the refinement.
As demonstrated later in \cref{sec:issue}, this expectation is wrong.
Immediately after a refinement with interpolation as transfer, the time error estimator of the implicit Euler method substantially overestimates the local error, and with adaptive step size control, this leads to a substantial step size reduction, caused solely by the perturbation introduced by the mesh change.
Such step size reductions after refinement with piecewise linear interpolation were observed in~\cite{BCJ98}.
The authors traced them back to the residual of the semi-discrete system on the new mesh and avoided them in one spatial dimension by cubic spline interpolation.

Adaptive methods for parabolic problems with meshes that change in time, controlled by \emph{a~posteriori} error estimates, were analyzed, for instance, in~\cite{EJ91,Pic98,LM06,SV08,KMSS12,GLW21}.
The elliptic reconstruction, introduced in~\cite{MN03} for the \emph{a~posteriori} analysis of spatially semi-discrete schemes and in~\cite{LM06} for fully discrete schemes, maps a finite element function to the exact solution of the elliptic problem with its discrete Laplacian as right-hand side.
The transfer proposed in this work maps the old discrete solution to the finite element approximation of its elliptic reconstruction in the new space, and thus, under refinement, it preserves the discrete Laplacian, whereas interpolation preserves the solution.
Its variant expresses the discrete Laplacian through the implicit Euler scheme, with the data of the differential equation instead of their projection.
It coincides with the elliptic transfer operator of~\cite{CGS21}, which leads to \emph{a~posteriori} error bounds of optimal order for Galerkin spaces that change arbitrarily between time steps.
For the Crank--Nicolson method, which also evaluates the discrete Laplacian at the previous time step, it was shown in~\cite{BKM12,BKM13} that the discrete Laplacian of the old solution on the new mesh causes oscillations after refinement, and that replacing it by the $L^2$ projection of the old discrete Laplacian avoids them, while the solution itself is still transferred unchanged.
It was also observed that their time error estimators jump at spatial refinement when the solution varies rapidly in space, and the study of nonstandard transfers was suggested~\cite[Sections~3.4 and~4]{BKM13}.
For the transfer proposed in this work, the discrete Laplacian of the transferred solution is the $L^2$ projection of the old one, so that no separate transfer of the discrete Laplacian is needed, as also noted for the elliptic transfer operator in~\cite{CGS21}.
To the best of our knowledge, apart from the formal one-dimensional analysis in~\cite{BCJ98}, the effect of the mesh transfer on the local time error estimator, which controls the step size, has not been analyzed.

Our contributions are as follows.
After an arbitrary mesh change, we derive an identity for the discrete time derivative of the first step on the new mesh that isolates the transfer residual,~i.e.,~the defect of the discrete equation on the new mesh for the transferred solution and time derivative.
Unless this residual vanishes, we prove two-sided bounds showing that the estimator overestimates the local error by a factor proportional to the inverse step size for small step sizes.
For interpolation and the $L^2$ projection, the residual after refinement is, up to the data oscillation, the component of the discrete Laplacian of the old solution on the new mesh that is orthogonal to the old space.
A transfer based on the elliptic reconstruction reduces the residual after refinement to the data oscillation, and its variant with the data of the differential equation makes it vanish for arbitrary mesh changes if the discrete time derivative is transferred by the $L^2$ projection.
In a numerical experiment with a single global refinement, we first confirm this overestimation for interpolation and its dependence on the step size.
Afterwards, we demonstrate in a space and time adaptive experiment the advantages of the nonstandard mesh transfers in comparison to interpolation.

The paper is organized as follows.
In \cref{sec:setting}, we introduce the notation, the implicit Euler method with its time error estimator and the mesh transfer, and demonstrate that the estimator overestimates the local error after a refinement with interpolation.
In \cref{sec:first-step}, we analyze the first time step after a mesh change.
In \cref{sec:transfers}, we determine the transfer residual of interpolation, the $L^2$ projection, two transfers based on the elliptic reconstruction, and a constrained variant that preserves the old nodal values.
In \cref{sec:algorithm}, we describe a space and time adaptive algorithm with step size control and spatial mesh adaptation.
Finally, \Cref{sec:numerics} contains numerical experiments with one refinement and with a moving bump computed with adaptivity in space and time, and \cref{sec:conclusions} concludes.

\section{Setting}\label{sec:setting}

In this section, we introduce the model problem, its discretization in space and time, the time error estimator and the mesh transfer.
Moreover, we collect some properties of the discrete operators under nested refinement, and demonstrate in an example how refinement affects the time error estimator.

\subsection{Model problem and notation}\label{sec:model}
Let $\Omega\subset\R^m$, $m\in\{1,2,3\}$, be a bounded convex polygonal domain, and consider the heat equation
\begin{equation}\label{eq:heat}
  \partial_t u - \Delta u = f \quad \text{in }\Omega\times(0,T),
  \qquad
  u=0 \quad \text{on }\partial\Omega\times(0,T),
  \qquad
  u(0)=u_0 \quad \text{in }\Omega .
\end{equation}
For $k\in\mathbb{N}$, we denote by $H^k(\Omega)$ the standard Sobolev spaces with the norms $\norm{\cdot}_{H^k(\Omega)}$, and by $H^1_0(\Omega)$ the subspace of functions in $H^1(\Omega)$ with vanishing trace on $\partial\Omega$.
We denote by $(\cdot,\cdot)$ and $\norm\cdot$ the inner product and the norm of $L^2(\Omega)$ and set $a(v,w):=(\nabla v,\nabla w)$.
For an operator $B: V \to W$ with $V, W \subset L^2(\Omega)$, we denote
by $\norm{B} = \sup_{0 \neq v \in V} \norm{Bv} / \norm{v}$ the operator norm induced
by the $L^2$ norm, and by $\Id$ the identity.

\subsection{Finite element spaces and discrete operators}\label{sec:fe}

For a shape-regular mesh $\T_X$ of $\Omega$, we denote by $X := \max_{K \in \T_X} \diam(K)$ the maximal mesh size and by $V_X \subset H_0^1(\Omega)$ the Lagrange finite element space of degree $p$.
We denote by $P_X:L^2(\Omega)\to V_X$ the $L^2$ projection, and define the discrete Laplacian $A_X:H^1_0(\Omega)\to V_X$ and the Ritz projection $R_X:H^1_0(\Omega)\to V_X$ by
\begin{align*}
  (A_Xv,w)=a(v,w),\qquad a(R_Xv,w)=a(v,w)\qquad\text{for all }w\in V_X.
\end{align*}
Moreover, we denote by $\N_X$ the set of Lagrange nodes of $V_X$, including those on $\partial\Omega$, and by $I_X$ the nodal interpolation, which maps $v\in C(\bar\Omega)$ with $v=0$ on $\partial\Omega$ to the function $I_Xv\in V_X$ with $(I_Xv)(z)=v(z)$ for all $z\in\N_X$.
Restricted to $V_X$, $A_X$ is self-adjoint and positive definite with respect to $(\cdot,\cdot)$.
Since $V_X$ is finite dimensional, it has an $L^2$-orthonormal basis of eigenfunctions, and we denote by $\lambda_X:=\max_{0\neq v\in V_X}\norm{\nabla v}^2/\norm{v}^2=\norm{A_X|_{V_X}}$ its largest eigenvalue.
The \emph{elliptic reconstruction}~\cite{MN03} of $v\in V_X$ is defined by $\Rec_Xv\in H^1_0(\Omega)$ with
\begin{align*}
    a(\Rec_Xv,\varphi)=(A_Xv,\varphi)\qquad\text{for all }\varphi\in H^1_0(\Omega).
\end{align*}

\subsection{Implicit Euler method and time error estimator}\label{sec:time}
Let $0 = t_0 < t_1 < \dots$ be a time grid with step sizes $\tau_n = t_n - t_{n-1}$, $n \ge 1$, and $f^n := f(t_n)$ for $n \ge 0$.
On a fixed mesh $\T_X$ and for an approximation $u_X^0\in V_X$ of $u_0$, the implicit Euler method computes $u_X^n \in V_X$, such that
\begin{equation}\label{eq:euler} 
  d_X^n+A_Xu_X^n=P_Xf^n,
  \qquad
  d_X^n:=\frac{u_X^n-u_X^{n-1}}{\tau_n},
\end{equation}
for $n \ge 1$, where $d_X^n$ is the discrete time derivative.
The scheme \eqref{eq:euler} is the finite element discretization of the semi-discretization of \eqref{eq:heat} in time by the implicit Euler method, which finds $U^n \in H_0^1(\Omega)$ with $U^0=u_0$ by
\begin{equation}\label{eq:semi-euler}
  U^n-\tau_n\Delta U^n=U^{n-1}+\tau_nf^n
\end{equation}
for $n \ge 1$.
Inserting the exact solution into \eqref{eq:semi-euler} and Taylor expanding $u(t_{n-1})$ about $t_n$ yields the local truncation error (LTE)
\begin{equation}\label{eq:euler-ex-est}
  u(t_n)-u(t_{n-1})-\tau_n \bigl(\Delta u(t_n)+f(t_n)\bigr)
  =-\frac{\tau_n^2}2\,\partial_{tt}u(t_n)+\OO(\tau_n^3).
\end{equation}
Approximating $\partial_{tt} u(t_n)$ with the second divided difference, we obtain for the norm of the leading term in \eqref{eq:euler-ex-est} the approximation
\begin{equation}\label{eq:est-euler}
  \eta^n:=\kappa_n\tau_n\norm{d_X^n-d_X^{n-1}},\qquad
  \kappa_n:=\frac{\tau_n}{\tau_n+\tau_{n-1}} .
\end{equation}
For $n=1$, we use $d_X^0:=P_Xf^0-A_Xu_X^0$,~i.e.,~the first equation in~\eqref{eq:euler} for $n=0$, as an approximation of $\partial_tu(t_0)$, and set $\tau_0:=0$.
To control the step size adaptively, we prescribe a tolerance $\tol \mathrel{>} 0$ and require $\eta^n \le \tol$.
A common choice for the next step size, see~\cite[Section~II.4]{HNW93}, is then given by
\begin{equation}\label{eq:controller-intro}
  \tau_{n+1}=\tau_n \left(\frac{\tol}{\eta^n}\right)^{1/2}.
\end{equation}

\subsection{Mesh change}\label{sec:mesh-change}
Suppose that the mesh $\T_H$ is replaced by $\T_h$ after step $N$, so that $\tau_N$ is the last step size on the old mesh and $\tau_{N+1}$ the first step size on the new mesh.
We denote by $u_H^n,d_H^n\in V_H$, $n\le N$, the solutions of~\eqref{eq:euler} on $\T_H$, and by $u_h^n,d_h^n\in V_h$, $n\ge N+2$, those on $\T_h$.
The first step on the new mesh determines $u_h^{N+1}$ from the transferred solution $\Pi u_H^N \in V_h$ through
\begin{equation}\label{eq:intro-transfer}
  d_h^{N+1}+A_hu_h^{N+1}=P_hf^{N+1},\qquad
  d_h^{N+1}:=\frac{u_h^{N+1}-\Pi u_H^N}{\tau_{N+1}},
\end{equation}
where $\Pi:V_H\to V_h$ is a suitable \emph{mesh transfer}.
We refer to $V_H$ as the old and $V_h$ as the new space.
Unless stated otherwise, no relation between the two is assumed.
We speak of \emph{nested refinement} if $\T_h$ is obtained by refining cells of $\T_H$.
For conforming Lagrange finite elements on affine simplicial meshes, this implies $V_H\subseteq V_h$.
This inclusion may fail for other discretizations, e.g., for nonconforming elements.
For nested refinement, we have the following properties of the introduced operators.

\begin{lemma}\label{lem:nesting}
  If $V_H\subset V_h$, then for all $v\in H^1_0(\Omega)$ and $v_H\in V_H$, it holds
  \begin{align*}
      P_HP_hv=P_Hv, \qquad P_HA_hv=A_Hv, \qquad I_hv_H=P_hv_H=R_hv_H=v_H .
  \end{align*}
\end{lemma}
\begin{proof}
  Let $w\in V_H$.
  Since $V_H\subset V_h$, we obtain
  \begin{align*}
    (P_HP_hv,w)&=(P_hv,w)=(v,w)=(P_Hv,w),\\
    (P_HA_hv,w)&=(A_hv,w)=a(v,w)=(A_Hv,w).
  \end{align*}
  Hence, $P_HP_hv-P_Hv\in V_H$ and $P_HA_hv-A_Hv\in V_H$ are orthogonal to $V_H$, and thus vanish.
  This shows the first two identities.
  For the remaining ones, we use $v_H\in V_H \subset V_h$.
  By the definitions of $P_h$ and $R_h$, we have $(P_hv_H-v_H,w)=0$ and $a(R_hv_H-v_H,w)=0$ for all $w\in V_h$.
  Choosing $w=P_hv_H-v_H\in V_h$ and $w=R_hv_H-v_H\in V_h$, respectively, yields $\norm{P_hv_H-v_H}=0$ and $\norm{\nabla(R_hv_H-v_H)}=0$.
  Since $\norm{\nabla\cdot}$ is a norm on $H^1_0(\Omega)$, we obtain $P_hv_H=R_hv_H=v_H$.
  Finally, $I_hv_H$ and $v_H$ are functions in $V_h$ with the same values at all nodes of $\T_h$, and thus coincide.
\end{proof}

\subsection{The time error estimator after a refinement}\label{sec:issue}
\begin{figure}[htb!]
  \centering
  \begin{tikzpicture}
  \definecolor{figink}{HTML}{52514e}
  \definecolor{figmuted}{HTML}{898781}
  \definecolor{figgrid}{HTML}{e1e0d9}
  \definecolor{figaxis}{HTML}{c3c2b7}
  \definecolor{figinterp}{HTML}{cc3311}%
  \pgfplotstableread[row sep=\\]{
    n est exact\\
    1 nan 2.49975e-09\\
    2 2.65525e-09 2.4995e-09\\
    3 2.65468e-09 2.49925e-09\\
    4 2.6541e-09 2.499e-09\\
    5 2.65353e-09 2.49875e-09\\
    6 2.65296e-09 2.4985e-09\\
    7 2.65239e-09 2.49825e-09\\
    8 2.65182e-09 2.498e-09\\
    9 2.65125e-09 2.49775e-09\\
    10 2.65068e-09 2.4975e-09\\
    11 0.000270506 2.49725e-09\\
    12 0.000148735 2.497e-09\\
    13 6.69687e-05 2.49675e-09\\
    14 3.01531e-05 2.4965e-09\\
    15 1.35766e-05 2.49625e-09\\
    16 6.11298e-06 2.496e-09\\
    17 2.75241e-06 2.49575e-09\\
    18 1.23929e-06 2.4955e-09\\
    19 5.57999e-07 2.49525e-09\\
    20 2.51243e-07 2.495e-09\\
    21 1.13124e-07 2.49476e-09\\
    22 5.09353e-08 2.49451e-09\\
    23 2.2935e-08 2.49426e-09\\
    24 1.0329e-08 2.49401e-09\\
    25 4.65603e-09 2.49376e-09\\
    26 2.10863e-09 2.49351e-09\\
    27 9.76991e-10 2.49326e-09\\
    28 4.98331e-10 2.49301e-09\\
    29 3.2698e-10 2.49276e-09\\
    30 2.82873e-10 2.49251e-09\\
    31 2.7633e-10 2.49226e-09\\
    32 2.7827e-10 2.49201e-09\\
    33 2.81945e-10 2.49176e-09\\
    34 2.85962e-10 2.49151e-09\\
    35 2.90041e-10 2.49127e-09\\
    36 2.94125e-10 2.49102e-09\\
    37 2.98203e-10 2.49077e-09\\
    38 3.02273e-10 2.49052e-09\\
    39 3.06335e-10 2.49027e-09\\
    40 3.10388e-10 2.49002e-09\\
  }\problemheatdata
  \begin{semilogyaxis}[
      scale only axis, width=6.1cm, height=4.2cm,
      xmin=0, xmax=40, xtick distance=5,
      ymin=1e-10, ymax=1e-3, ytick={1e-10,1e-9,1e-8,1e-7,1e-6,1e-5,1e-4,1e-3},
      axis lines*=left, axis line style={-, figaxis, line width=0.6pt},
      tick align=outside, tick style={figaxis},
      grid=major, major grid style={figgrid, line width=0.4pt},
      tick label style={font=\scriptsize, text=figink},
      label style={font=\footnotesize, text=figink},
      xlabel={Step $n$},
      legend style={at={(0.5,-0.3)}, anchor=north, draw=none, fill=none,
        font=\footnotesize, legend columns=2,
        /tikz/every odd column/.append style={column sep=0.2em},
        /tikz/every even column/.append style={column sep=1.5em}},
      legend cell align=left,
    ]
    \draw[black, line width=1pt]
      ({axis cs:10.5,1}|-{rel axis cs:0,0}) -- ({axis cs:10.5,1}|-{rel axis cs:0,1})
      node[anchor=north east, xshift=-3pt, yshift=-3pt, inner sep=0pt, align=left,
        font=\scriptsize, text=black] {refined};
    \addplot[figink, line width=0.9pt, dash pattern=on 3.5pt off 1.75pt]
      table[x=n, y=exact] {\problemheatdata};
    \addlegendentry{local truncation error}
    \addplot[figinterp, line width=0.8pt, line join=round,
        mark=*, mark size=1.2pt, mark options={solid}]
      table[x=n, y=est] {\problemheatdata};
    \addlegendentry{estimator}
  \end{semilogyaxis}
\end{tikzpicture}
  \caption{
  Local error estimator \eqref{eq:est-euler} and LTE \eqref{eq:euler-ex-est} for the implicit Euler discretization of the heat equation with step size $\tau=10^{-4}$ and finite element degree $p=1$.
  After a single global mesh refinement, the solution is transferred to the fine mesh by interpolation.
  }
  \label{fig:problem}
\end{figure}
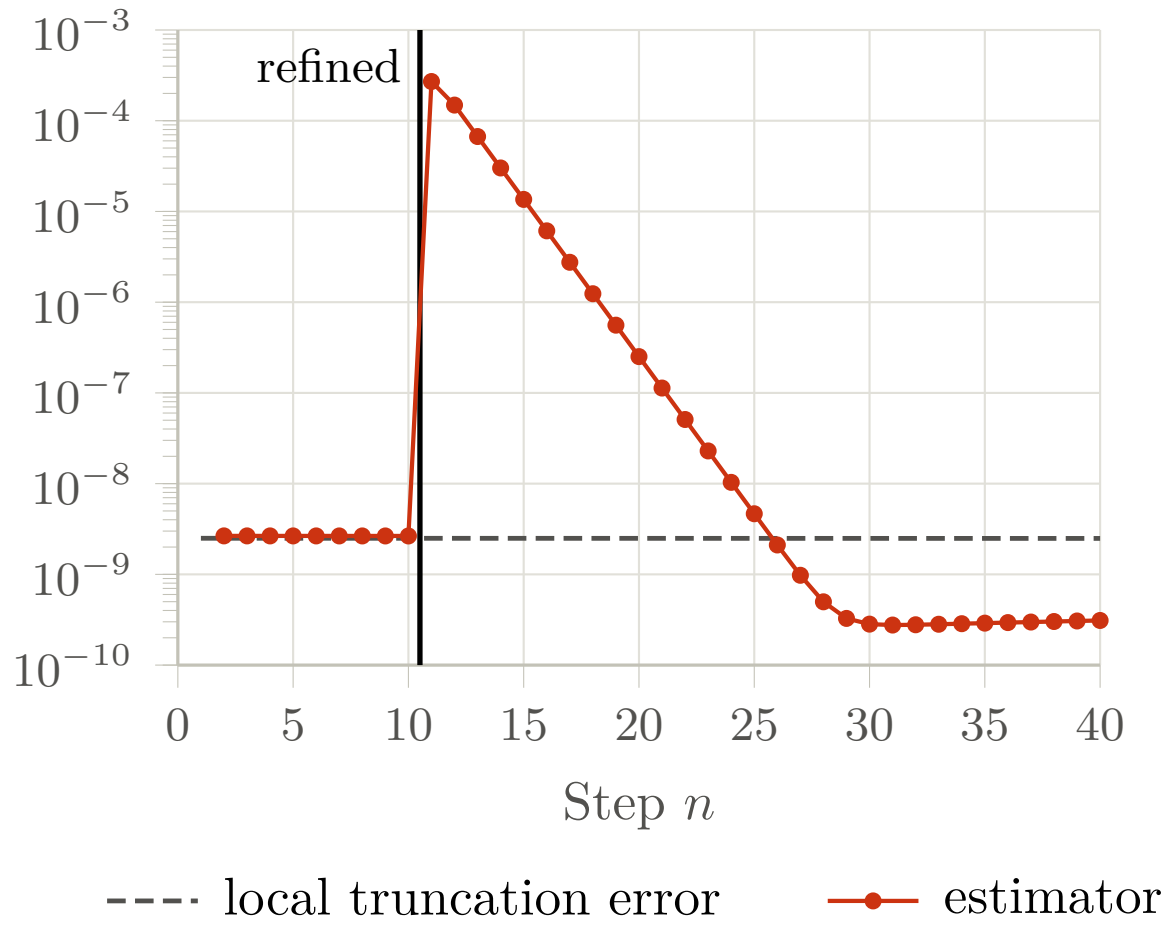

Under nested refinement, the nodal interpolation $I_h$ represents the solution exactly,~i.e.,~$\Pi u_H^N=u_H^N$ for $\Pi=I_h$ by \cref{lem:nesting}.
One may therefore expect the estimator \eqref{eq:est-euler} to remain unaffected by the refinement.
However, \cref{fig:problem} demonstrates that this expectation is wrong.
We consider $\Omega = (0, 1)^2$ and the smooth solution $u(x,t)=e^{-t}\sin(\pi x_1)\sin(\pi x_2)$ with the right-hand side $f(x,t)=(2\pi^2-1)u(x,t)$ for $x=(x_1,x_2)\in\Omega$, and $p=1$.
The uniform mesh with $H=\frac1{16}$ is refined globally after step $N=10$ to $h=\frac1{32}$, and the solution is transferred by interpolation.
In \cref{fig:problem}, we compare the local error estimator \eqref{eq:est-euler} with the LTE \eqref{eq:euler-ex-est} for the constant step size $\tau=10^{-4}$.
Before refinement, the estimator agrees with the LTE.
Immediately after refinement, it substantially overestimates the LTE.
With adaptive step size control as in \eqref{eq:controller-intro}, this increase in $\eta^n$ would lead to a substantial step size reduction, caused solely by the perturbation introduced by the mesh change.

\section{The first step after a mesh change}\label{sec:first-step}

In this section, we analyze the first implicit Euler step~\eqref{eq:intro-transfer} after a spatial mesh change, allowing general transfers of the discrete solution and its discrete time derivative.
In the estimator \eqref{eq:est-euler}, we compare $d_h^{N+1}$ with the old discrete time derivative $d_H^N$, transferred by a possibly different mesh transfer $\Pi_d:V_H\to V_h$,~i.e.,
\begin{align}\label{EQ:est-transfer}
    \eta^{N+1}=\kappa_{N+1}\tau_{N+1}\norm{d_h^{N+1}-\Pi_d d_H^N}.
\end{align}
Recall that the discrete time derivatives are defined in {\eqref{eq:euler} and~\eqref{eq:intro-transfer} by
\begin{align*}
    d_H^N = \frac{u_H^{N}-u_H^{N-1}}{\tau_{N}}, \qquad d_h^{N+1}=\frac{u_h^{N+1}-\Pi u_H^N}{\tau_{N+1}}.
\end{align*}
The solution $u_h^{N+1}$ depends only on $\Pi u_H^N$, while the estimator also depends on $\Pi_d d_H^N$.
To describe the effect of the mesh change from step $N$ to $N+1$, we introduce the transfer residual as follows.

\begin{definition}
     Let $\Pi,\Pi_d:V_H\to V_h$ be mesh transfers. 
     The \emph{transfer residual} is defined by
    \begin{equation}\label{eq:residual}
        r_N:=P_hf^N-\Pi_d d_H^N-A_h\Pi u_H^N\in V_h .
    \end{equation}
\end{definition}

This residual measures the defect of the discrete equation \eqref{eq:euler} on the new mesh at $t_N$ for the transferred discrete derivative and solution.
In the absence of a mesh change, \eqref{eq:euler} yields directly $r_N=0$.
If $\Pi_d=\Pi$ is linear, then $\Pi_d d_H^N=(\Pi u_H^N-\Pi u_H^{N-1})/\tau_N$.
In this case, $r_N$ is the defect of the transferred pair $(\Pi u_H^{N-1},\Pi u_H^N)$.
The transfer residual gives the following identity for the discrete time derivative after a mesh change.

\begin{lemma}\label{lem:first-step}
  Let $u_H^{N-1}, u_H^N\in V_H$ satisfy \eqref{eq:euler} on $\T_H$ for $n=N$, let $\Pi,\Pi_d:V_H\to V_h$ be mesh transfers, let $u_h^{N+1}\in V_h$ satisfy \eqref{eq:intro-transfer}, and let $r_N$ be the transfer residual given by \eqref{eq:residual}.
  Set $B_h^n = I + \tau_n A_h$.
  Then
  \begin{equation}\label{eq:first-step}
    B_h^{N+1} d_h^{N+1}=\Pi_d d_H^N+P_h(f^{N+1}-f^N)+r_N .
  \end{equation}
\end{lemma}
\begin{proof}
  By \eqref{eq:intro-transfer} and $u_h^{N+1}=\Pi u_H^N+\tau_{N+1}d_h^{N+1}$,
  we have $B_h^{N+1}d_h^{N+1}=P_hf^{N+1}-A_h\Pi u_H^N$. Inserting $A_h\Pi u_H^N=P_hf^N
  -\Pi_d d_H^N-r_N$ from \eqref{eq:residual} yields \eqref{eq:first-step}.
\end{proof}

Next, using the identity \eqref{eq:first-step}, we separate the contributions to the error estimator \eqref{EQ:est-transfer} as follows.

\begin{corollary}\label{thm:splitting}
  Under the assumptions of \cref{lem:first-step}, let $S_h^{N+1} := (B_h^{N+1})^{-1} = (\Id+\tau_{N+1} A_h)^{-1}$, and define
  \begin{equation}\label{eq:split-terms}
    \delta_d:=A_h\Pi_d d_H^N-P_hA_Hd_H^N,\qquad
    \Phi:=S_h^{N+1}\bigl[P_h(f^{N+1}-f^N)-\tau_{N+1}P_hA_Hd_H^N\bigr].
  \end{equation}
  Then, it holds
  \begin{equation}\label{eq:splitting}
    d_h^{N+1}-\Pi_d d_H^N=\Phi-\tau_{N+1}S_h^{N+1}\,\delta_d+S_h^{N+1}\,r_N .
  \end{equation}
  Without a mesh change, $\delta_d=0$ and $r_N=0$.
\end{corollary}
\begin{proof}
  Subtracting $B_h^{N+1}\Pi_d d_H^N$ from \eqref{eq:first-step} and inserting
  $A_h\Pi_d d_H^N=P_hA_Hd_H^N+\delta_d$ gives
  \begin{align*}
    B_h^{N+1}(d_h^{N+1}-\Pi_d d_H^N)=P_h(f^{N+1}-f^N)-\tau_{N+1}P_hA_Hd_H^N
    -\tau_{N+1}\delta_d+r_N .
  \end{align*}
  Applying $S_h^{N+1}$ yields \eqref{eq:splitting}.
\end{proof}

We define the second-order difference quotient across a mesh transfer by 
\begin{align*}
  D_h^{N+1}:=\frac{d_h^{N+1}-\Pi_d d_H^N}{\tau_{N+1}},
\end{align*}
which approximates $\partial_{tt}u(t_{N+1})$ on a fixed spatial mesh.
The estimator~\eqref{EQ:est-transfer} satisfies $\eta^{N+1} = \kappa_{N+1} \tau_{N+1}^2 \norm{D_h^{N+1}}$.
By~\eqref{eq:splitting}, the second-order difference quotient satisfies
\begin{align}\label{eq:dn1}
    D_h^{N+1} = \tau_{N+1}^{-1} \Phi - S_h^{N+1}\,\delta_d + \tau_{N+1}^{-1} S_h^{N+1}\,r_N.
\end{align}
Unless $r_N = 0$, the transfer residual contributes a term of the order of $\tau_{N+1}^{-1}$ for small $\tau_{N+1}$,
as stated in the following theorem. 

\begin{theorem}\label{thm:overestimation}
  Under the assumptions of \cref{thm:splitting}, let
  \begin{align*}
    M^{N+1}:=\sup_{t_N\le t\le t_{N+1}}\norm{\partial_tf(t)}+\norm{A_Hd_H^N}+\norm{\delta_d}.
  \end{align*}
  Then, for all $\tau_{N+1}>0$, it holds
  \begin{equation}\label{eq:D-bounds}
    \frac{\norm{r_N}}{(1+\tau_{N+1}\lambda_h)\tau_{N+1}}-M^{N+1}\le\norm{D_h^{N+1}}
    \le\frac{\norm{r_N}}{\tau_{N+1}}+M^{N+1} .
  \end{equation}
  If, in addition, $\tau_{N+1}\lambda_h\le1$ and $4\tau_{N+1}M^{N+1}\le\norm{r_N}$,
  then
  \begin{equation}\label{eq:D-tau1}
    \frac{\norm{r_N}}{4}\,\tau_{N+1}^{-1}\le\norm{D_h^{N+1}}
    \le\frac{5\norm{r_N}}{4}\,\tau_{N+1}^{-1}.
  \end{equation}
\end{theorem}
\begin{proof}
  Since $A_h$ is self-adjoint and positive definite on $V_h$, there is an
  $L^2$-orthonormal basis $(\varphi_k)_k$ of $V_h$ with $A_h\varphi_k=\lambda_k
  \varphi_k$ and $0<\lambda_k\le\lambda_h$, and thus $S_h^{N+1}\varphi_k=(1+\tau_{N+1}
  \lambda_k)^{-1}\varphi_k$. 
  For $v=\sum_kc_k\varphi_k\in V_h$, this
  yields, using Parseval's identity, $\norm{S_h^{N+1}v}^2=\sum_k(1+\tau_{N+1}\lambda_k)^{-2}c_k^2$.
  Hence, we obtain
  \begin{align*}
    \frac{\norm v}{1+\tau_{N+1}\lambda_h}\le\norm{S_h^{N+1}v}\le\norm v
    \qquad\text{for all }v\in V_h.
  \end{align*}
  Thus the bound holds for $v = r_N \in V_h$, and we have $\norm{S_h^{N+1}} \le 1$.
  With $\norm{P_h}\le1$ and $f^{N+1}-f^N=\int_{t_N}^{t_{N+1}}\partial_tf\,dt$, we obtain
  \begin{align*}
    \norm{\Phi-\tau_{N+1}S_h^{N+1}\delta_d}\le\norm{f^{N+1}-f^N}+\tau_{N+1}\norm{A_Hd_H^N}
    +\tau_{N+1}\norm{\delta_d}\le\tau_{N+1}M^{N+1} .
  \end{align*}
  Combined with \eqref{eq:dn1}, this concludes \eqref{eq:D-bounds}. 
  Using the additional conditions, we have $1+\tau_{N+1}\lambda_h\le2$ and $M^{N+1}\le\norm{r_N}/(4\tau_{N+1})$.
  This yields \eqref{eq:D-tau1}. 
\end{proof}

The proof of \cref{thm:overestimation} shows that the first two terms on the right-hand side of \eqref{eq:dn1} satisfy $\norm{\tau_{N+1}^{-1}\Phi-S_h^{N+1}\,\delta_d}\le M^{N+1}$, which remains bounded as $\tau_{N+1}\to0$. 
The residual $r_N$ does not depend on $\tau_{N+1}$. 
Hence, if $r_N\neq0$, its contribution dominates $D_h^{N+1}$ for small $\tau_{N+1}$.
If $r_N\perp V_H$, as shown later under nested refinement for two of the transfers in~\cref{lem:projection,prop:ER}, the following lemma shows that the factor $S_h^{N+1}$ in \eqref{eq:dn1} reduces this contribution by $H^2$ for fixed $\tau_{N+1}$ but still grows as $\tau_{N+1}\to0$.

\begin{lemma}\label{lem:damping}
  Let $V_H\subset V_h$, let $\Omega$ be convex and $\T_H$, $\T_h$ shape regular.
  Then, there exists a constant $C_D>0$, depending only on $\Omega$ and the 
  shape regularity, such that
  for all $\delta\in V_h$ with $\delta\perp V_H$ and all $\tau>0$ it holds
  \begin{align*}
    \norm{(\Id+\tau A_h)^{-1}\delta}\le\min\bigl\{1,C_DH^2/\tau\bigr\}\norm\delta .
  \end{align*}
\end{lemma}
\begin{proof}
  The bound with constant $1$ holds since $\norm{(\Id+\tau A_h)^{-1}}\le1$ on $V_h$,~cf.~the proof of \cref{thm:overestimation}. 
  Let $\varphi:=A_h^{-1}\delta\in V_h$ and let $\psi\in H^2(\Omega)\cap H^1_0(\Omega)$ solve $-\Delta\psi=\varphi$, so that   $\norm{\psi}_{H^2(\Omega)}\le C_1 \norm\varphi$~\cite[Theorem~3.2.1.2]{Gri85}. 
  Using $\norm{R_X\psi-\psi}\le C_2 X^2\norm\psi_{H^2(\Omega)}$, $X\in\{H,h\}$~\cite[Chapter~1]{Tho06}, we have $\norm{R_h \psi - R_H \psi} \le \norm{R_h \psi - \psi} + \norm{R_H \psi - \psi} \le C_2 (h^2 + H^2) \norm{\psi}_{H^2(\Omega)}$. 
  With $\delta\perp V_H$ and $h \le H$, we obtain
  \begin{align*}
    \norm\varphi^2=a(\psi,\varphi)=a(R_h\psi,\varphi)=(\delta,R_h\psi)
    =(\delta,R_h\psi-R_H\psi)\le C_D H^2\norm\delta\norm\varphi , 
  \end{align*}
  where $C_D = 2 C_1 C_2$.
  Further, we have $\norm{(\Id+\tau A_h)^{-1}A_h}\le\tau^{-1}$, since its eigenvalues are $\lambda_k/(1+\tau\lambda_k)<\tau^{-1}$.
  Using the above inequalities, we conclude
  \begin{align*}
      \norm{(\Id+\tau A_h)^{-1}\delta}=\norm{(\Id+\tau A_h)^{-1}A_h\varphi}\le \tau^{-1}\norm\varphi\le C_DH^2\tau^{-1}\norm\delta .
  \end{align*}
\end{proof}

\section{Mesh transfer operators}\label{sec:transfers}

In this section, we determine the transfer residual~\eqref{eq:residual} for interpolation and the $L^2$ and Ritz projections, which preserve the old solution under nested refinement, and for three transfers based on the elliptic reconstruction.

\subsection{Interpolation and \texorpdfstring{$L^2$}{L2} projection}\label{sec:interpolation-trans}
A common approach is the nodal interpolation operator $I_h: V_H \to V_h$ with $\Pi = \Pi_d = I_h$.
Under nested refinement $V_H \subset V_h$, this transfer satisfies $\Pi v = v$ for all $v \in V_H$.
The same property holds for the $L^2$ projection $P_h$ and the Ritz projection $R_h$,~cf.~\cref{lem:nesting}.
The $L^2$ projection arises implicitly if the old solution enters the scheme only through the $L^2$ inner product with the test functions of the new space, as in the discontinuous Galerkin time discretization of~\cite{EJ91} and in the implicit Euler method of~\cite{LM06,KMSS12}.
The Ritz projection is used as transfer in~\cite{DGMW26}.
For these transfers, the following lemma characterizes the corresponding transfer residual for nested refinement.

\begin{lemma}\label{lem:projection}
  Let $V_H\subset V_h$ with shape-regular meshes $\T_H$ and $\T_h$, and let $\Pi,\Pi_d:V_H\to V_h$ be mesh transfers with $\Pi v=\Pi_d v=v$ for all $v\in V_H$.
  Then, the following holds.
  \begin{enumerate}[label=(\roman*)]
    \item Let $u_H^N,d_H^N\in V_H$ satisfy \eqref{eq:euler} on $\T_H$ for $n=N$.
      Then, the transfer residual is given by
      \begin{equation}\label{eq:residual-projection}
        r_N=(P_h-P_H)f^N-\delta_N,
        \qquad
        \delta_N:=(\Id-P_H)A_hu_H^N.
      \end{equation}
      Moreover, $r_N\perp V_H$ and $\delta_N\perp V_H$.

    \item Let $\T_h$ and $\T_H$ be quasi-uniform.
      For every $v\in V_H$, let $\omega:=\Rec_Hv$.
      Then,
      \begin{equation}\label{eq:delta-bound}
        \frac{\norm{\nabla(R_h-R_H)\omega}}{C_AH}
        \le\norm{(\Id-P_H)A_hv}
        \le\frac{C_{\mathrm{inv}}}{h}
        \norm{\nabla(R_h-R_H)\omega},
      \end{equation}
      where $C_A$ and $C_{\mathrm{inv}}$ are the constants in
      \begin{align*}
        \norm{\varphi-P_H\varphi}
        \le C_AH\norm{\nabla\varphi},
        \qquad \varphi\in H_0^1(\Omega),
      \end{align*}
      and
      \begin{align*}
        \norm{A_h\varphi}
        \le C_{\mathrm{inv}}h^{-1}\norm{\nabla\varphi},
        \qquad \varphi\in V_h.
      \end{align*}

    \item Under the assumptions of (ii), let $H\le ch$ for some $c>0$. 
    If, in addition, the saturation assumption
      \begin{align*}
        \norm{\nabla(\omega-R_h\omega)}
        \le\beta\norm{\nabla(\omega-R_H\omega)}
      \end{align*}
      holds with some $0\le\beta<1$, then, with $E_H(\omega):=\min_{w\in V_H}\norm{\nabla(\omega-w)}$,
      \begin{equation}\label{eq:delta-order}
        \frac{(1-\beta^2)^{1/2}}{C_A}
        \frac{E_H(\omega)}{H}
        \le\norm{(\Id-P_H)A_hv}
        \le C_{\mathrm{inv}}c\frac{E_H(\omega)}{H}.
      \end{equation}
  \end{enumerate}
\end{lemma}

\begin{proof}
  We start with (i).
  By \eqref{eq:euler} and \cref{lem:nesting}, we have
  \begin{align*}
    r_N
    &=P_hf^N-d_H^N-A_hu_H^N =P_hf^N-P_Hf^N+A_Hu_H^N-A_hu_H^N\\
    &=(P_h-P_H)f^N-(\Id-P_H)A_hu_H^N.
  \end{align*}
  Since $(\Id-P_H)g\perp V_H$ for every $g\in L^2(\Omega)$,
  we obtain $\delta_N\perp V_H$.
  Moreover, $(P_h-P_H)f^N=(\Id-P_H)P_hf^N$ by
  \cref{lem:nesting}, and hence $r_N\perp V_H$.

  We continue with (ii).
  Let $v\in V_H$ and $\omega=\Rec_Hv$.
  By the definition of the elliptic reconstruction,
  \begin{align*}
    a(\omega,\varphi)=(A_Hv,\varphi)
    \qquad\text{for all }\varphi\in H_0^1(\Omega).
  \end{align*}
  In particular, $R_H\omega=v$.
  For every $\varphi\in V_h$, we have
  \begin{align*}
    a(R_h\omega,\varphi)
    =a(\omega,\varphi)
    =(A_Hv,\varphi),
  \end{align*}
  and thus $A_hR_h\omega=A_Hv=P_HA_hv$ by
  \cref{lem:nesting}.
  Consequently,
  \begin{align*}
    (\Id-P_H)A_hv
    =A_h(v-R_h\omega)
    =A_h(R_H-R_h)\omega.
  \end{align*}
  Let $\delta:=(\Id-P_H)A_hv$ and
  $\varphi:=(R_H-R_h)\omega$.
  Since $\delta\perp V_H$, we obtain
  \begin{align*}
    \norm{\nabla\varphi}^2
    &=(\delta,\varphi)
    =(\delta,\varphi-P_H\varphi)\\
    &\le C_AH\norm\delta\norm{\nabla\varphi}.
  \end{align*}
  This gives the lower bound in \eqref{eq:delta-bound}.
  The upper bound follows from the inverse estimate
  \begin{align*}
    \norm\delta
    =\norm{A_h\varphi}
    \le C_{\mathrm{inv}}h^{-1}\norm{\nabla\varphi}.
  \end{align*}

  Finally, we prove (iii).
  By Galerkin orthogonality and $V_H\subset V_h$,
  \begin{align*}
    \norm{\nabla(\omega-R_H\omega)}^2
    =\norm{\nabla(\omega-R_h\omega)}^2
    +\norm{\nabla(R_h-R_H)\omega}^2.
  \end{align*}
  Since
  $\norm{\nabla(\omega-R_H\omega)}=E_H(\omega)$,
  the saturation assumption yields
  \begin{align*}
    (1-\beta^2)^{1/2}E_H(\omega)
    \le\norm{\nabla(R_h-R_H)\omega}
    \le E_H(\omega).
  \end{align*}
  Combining this with \eqref{eq:delta-bound} and
  $h^{-1}\le cH^{-1}$ gives \eqref{eq:delta-order}.
\end{proof}

Taking $v=u_H^N$ in \eqref{eq:delta-order} shows that the $L^2$ norm of the defect $\delta_N$ in the residual $r_N$ is, up to constants, bounded from above and from below by the best approximation error $E_H(\omega_N)$ in $V_H$ of the elliptic reconstruction $\omega_N:=\Rec_Hu_H^N$, divided by $H$.
For sufficiently smooth solutions and finite element degree $p$, the best approximation error satisfies $E_H(\omega_N)\in\OO(H^p)$.
Conversely, on quasi-uniform meshes, the best approximation error of a smooth function that is not a polynomial of degree $p$ is bounded from below, up to some constant, by $H^p$~\cite{LXX14}. 
Thus, the transfer residual for interpolation generally does not converge with respect to the spatial mesh size for $p=1$, and for $p=2$ only linearly.

The saturation assumption in \cref{lem:projection}~(iii) is known in the analysis of hierarchical \emph{a~posteriori} error estimators~\cite{BS93} and~\cite[Chapter~5]{AO00}.
It may fail for oscillatory data, but holds if the data oscillation is small compared to the error, as shown in~\cite{DN02} for the enrichment of linear by quadratic elements.
In \cref{lem:projection}~(iii), it is only needed for the lower bound, since the upper bound follows from \eqref{eq:delta-bound} and Galerkin orthogonality.

\subsection{Elliptic reconstruction transfer}\label{sec:ell-rec-trans}

Next, we consider two transfers based on the elliptic reconstruction, which remove the defect $\delta_N$ from the transfer residual.
For the discrete time derivative transfer, we choose $\Pi_d=P_h$.

\begin{definition}\label{def:ER}
  The elliptic reconstruction transfer (ERT) $\PiER:V_H\to V_h$ is defined by
  \begin{align*}
    a(\PiER v,w)=(A_Hv,w)
    \qquad\text{for all }w\in V_h.
  \end{align*}
\end{definition}

The definition of the ERT gives $A_h\PiER v=P_h A_Hv$ for all $v\in V_H$, and thus $A_h\PiER v=A_Hv$ if $V_H\subset V_h$.

\begin{proposition}\label{prop:ER}
  Let $\Pi = \PiER$ be the ERT and $\Pi_d = P_h$.
  Assume $V_H \subset V_h$ and shape regularity of $\T_H, \T_h$, let $u_H^N,d_H^N\in V_H$ satisfy \eqref{eq:euler} on $\T_H$ for $n=N$.
  Then, the transfer residual \eqref{eq:residual} satisfies
  \begin{align}\label{eq:resdiual-ert}
    r_N=P_h(\Id-P_H)f^N=(P_h-P_H)f^N.
  \end{align}
  In particular, $r_N\perp V_H$.
  If $f^N\in H^{p+1}(\Omega)\cap H_0^1(\Omega)$, then
  \begin{align*}
    \norm{r_N}
    \le C H^{p+1}\norm{f^N}_{H^{p+1}(\Omega)},
  \end{align*}
  where $C>0$ depends only on the shape regularity and the polynomial degree.
\end{proposition}

\begin{proof}
  By \eqref{eq:euler} and the definition of the ERT, we have
  \begin{align*}
    r_N &=P_hf^N-P_hd_H^N-A_h\PiER u_H^N =P_h(f^N-d_H^N-A_Hu_H^N)\\
    &=P_h(f^N-P_Hf^N)=(P_h-P_H)f^N,
  \end{align*}
  where we used $P_hP_H=P_H$ from \cref{lem:nesting}.
  Since $(P_h-P_H)f^N=(\Id-P_H)P_hf^N$, we obtain $r_N\perp V_H$.
  Moreover, using $\norm{P_h}\le1$ and the best-approximation property of $P_H$, we obtain
  \begin{align*}
    \norm{r_N}
    &\le\norm{(\Id-P_H)f^N}=\inf_{v\in V_H}\norm{f^N-v}\le C H^{p+1}\norm{f^N}_{H^{p+1}(\Omega)},
  \end{align*}
  where the last inequality follows from the interpolation error estimate for $v=I_Hf^N$ in~\cite[Theorem~4.4.20]{BS08}.
\end{proof}

\begin{remark}
  For nested refinement, comparing \eqref{eq:residual-projection}
  with \eqref{eq:resdiual-ert}, the corresponding transfer residuals satisfy
  \begin{align}\label{eq:comp-interpol-ert}
    r_N^{\mathrm{ERT}}
    =r_N^{\mathrm{int}}+\delta_N,
    \qquad
    \delta_N=(\Id-P_H)A_hu_H^N.
  \end{align}
  Hence, the ERT eliminates the discrete Laplacian mismatch from the transfer residual, leaving only the data oscillation $(P_h-P_H)f^N$.
\end{remark}

The transfer residual of the ERT arises since the old discrete equation contains $P_H f^N$ instead of $P_h f^N$.
Next, we consider a variant that uses the data of the differential equation instead of the old discrete operator.

\begin{definition}\label{def:ERf}
  The elliptic reconstruction transfer with data $f$ (ERT-f)
  $\PiERf u_H^N:=\hat u_h^N\in V_h$ is defined by
  \begin{align*}
    a(\hat u_h^N,v)=(f^N-d_H^N,v)
    \qquad\text{for all }v\in V_h.
  \end{align*}
\end{definition}

The transfers ERT and ERT-f coincide if $f^N\in V_H$ or if $V_h\subset V_H$,~i.e.,~for pure coarsening.
By the definition of the ERT-f, we have immediately $A_h\hat u_h^N=P_h(f^N-d_H^N)$. 
The ERT-f coincides with the elliptic transfer of~\cite[Definition~3.19]{CGS21}, which is applied there to moving mesh and virtual element methods.
In contrast to \cref{prop:ER}, the following proposition requires no relation between the old and the new space, and thus applies to arbitrary mesh changes.

\begin{proposition}\label{thm:ERf}
  Let $\Pi = \PiERf$ be the ERT-f, $\hat u_h^N:=\PiERf u_H^N$ and $\Pi_d = P_h$.
  Then, the transfer residual \eqref{eq:residual} vanishes,~i.e.,~$r_N=0$.
\end{proposition}
\begin{proof}
  By the definition of the transfer residual \eqref{eq:residual} and the ERT-f, we have
  \begin{align*}
    r_N &=P_hf^N-P_hd_H^N-A_h\hat u_h^N=P_h(f^N-d_H^N)-P_h(f^N-d_H^N)=0.
  \end{align*}
\end{proof}

\begin{figure}[htb!]
  \centering
  \begin{tikzpicture}
  \definecolor{figink}{HTML}{52514e}
  \definecolor{figgrid}{HTML}{e1e0d9}
  \definecolor{figaxis}{HTML}{c3c2b7}
  \definecolor{figinterp}{HTML}{cc3311}
  \definecolor{figert}{HTML}{cc79a7}
  \definecolor{figertf}{HTML}{1c5cab}
  \pgfplotstableread[row sep=\\]{
    x ui di ue de uf df\\
    0 0 0 0 0 0 0\\
    0.125 0.359239 5.68582 0.381947 -0.29531 0.382261 -0.377882\\
    0.25 0.695703 -11.4041 0.70623 -0.877152 0.706375 -0.731822\\
    0.375 0.86728 13.7268 0.922103 -0.712941 0.922859 -0.912288\\
    0.5 0.983872 -16.1279 0.99876 -1.24048 0.998965 -1.03495\\
    0.625 0.86728 13.7268 0.922103 -0.712941 0.922859 -0.912288\\
    0.75 0.695703 -11.4041 0.70623 -0.877152 0.706375 -0.731822\\
    0.875 0.359239 5.68582 0.381947 -0.29531 0.382261 -0.377882\\
    1 0 0 0 0 0 0\\
  }\exampleonedata
  \pgfplotsset{
    exampleone/.style={
      scale only axis, width=3.8cm, height=2.4cm,
      xmin=0, xmax=1, xtick={0,0.25,0.5,0.75,1},
      axis lines*=left, axis line style={-, figaxis, line width=0.6pt},
      tick align=outside, tick style={figaxis},
      grid=major, major grid style={figgrid, line width=0.4pt},
      tick label style={font=\scriptsize, text=figink},
      label style={font=\footnotesize, text=figink},
      title style={font=\footnotesize, text=figink},
      ylabel style={at={(axis description cs:-0.22,0.5)}, anchor=south},
    },
    exampleonetop/.style={exampleone, ymin=0, ymax=1.08, xticklabels={}},
    exampleonebottom/.style={exampleone,
      xticklabels={$0$,$\frac14$,$\frac12$,$\frac34$,$1$},
      xlabel={$x$}, xlabel style={yshift=0.6ex}},
    exampleoneexact/.style={figink, line width=0.9pt, dash pattern=on 3.5pt off 1.75pt,
      domain=0:1, samples=121},
    exampleoneline/.style={line width=0.8pt, line join=round},
    exampleoneinterp/.style={exampleoneline, figinterp, mark=*, mark size=1.8pt,
      mark options={solid, fill=figinterp, draw=white, line width=0.3pt}},
    exampleoneert/.style={exampleoneline, figert, mark=square*, mark size=2pt,
      mark options={solid, fill=figert, draw=white, line width=0.3pt}},
    exampleoneertf/.style={exampleoneline, figertf, mark=diamond*, mark size=2.8pt,
      mark options={solid, fill=figertf, draw=white, line width=0.3pt}},
  }
  \begin{axis}[exampleonetop, name=ui, ylabel={$u_h^{N+1}$ and $u$}]
    \addplot[exampleoneexact] {exp(-0.001)*sin(180*x)};
    \addplot[exampleoneinterp] table[x=x, y=ui] {\exampleonedata};
  \end{axis}
  \begin{axis}[exampleonetop, name=ue, at={(ui.east)}, anchor=west, xshift=1.05cm]
    \addplot[exampleoneexact] {exp(-0.001)*sin(180*x)};
    \addplot[exampleoneert] table[x=x, y=ue] {\exampleonedata};
  \end{axis}
  \begin{axis}[exampleonetop, name=uf, at={(ue.east)}, anchor=west, xshift=1.05cm]
    \addplot[exampleoneexact] {exp(-0.001)*sin(180*x)};
    \addplot[exampleoneertf] table[x=x, y=uf] {\exampleonedata};
  \end{axis}
  \begin{axis}[exampleonebottom, name=di, at={(ui.south)}, anchor=north, yshift=-0.35cm,
      ymin=-18.5, ymax=16.1, ylabel={$d_h^{N+1}$ and $\partial_tu$}]
    \addplot[exampleoneexact] {-exp(-0.001)*sin(180*x)};
    \addplot[exampleoneinterp] table[x=x, y=di] {\exampleonedata};
  \end{axis}
  \begin{axis}[exampleonebottom, name=de, at={(di.east)}, anchor=west, xshift=1.05cm,
      ymin=-1.34, ymax=0.0992,
      legend style={at={(0.5,-0.52)}, anchor=north, draw=none, fill=none,
        font=\footnotesize, legend columns=4,
        /tikz/every odd column/.append style={column sep=0.2em},
        /tikz/every even column/.append style={column sep=1.2em}},
      legend cell align=left]
    \addlegendimage{exampleoneexact}\addlegendentry{exact}
    \addlegendimage{exampleoneinterp}\addlegendentry{interpolation}
    \addlegendimage{exampleoneert}\addlegendentry{ERT}
    \addlegendimage{exampleoneertf}\addlegendentry{ERT-f}
    \addplot[exampleoneexact, forget plot] {-exp(-0.001)*sin(180*x)};
    \addplot[exampleoneert, forget plot] table[x=x, y=de] {\exampleonedata};
  \end{axis}
  \begin{axis}[exampleonebottom, name=df, at={(de.east)}, anchor=west, xshift=1.05cm,
      ymin=-1.34, ymax=0.0992]
    \addplot[exampleoneexact] {-exp(-0.001)*sin(180*x)};
    \addplot[exampleoneertf] table[x=x, y=df] {\exampleonedata};
  \end{axis}
\end{tikzpicture}
  \caption{
  One step of the implicit Euler method with $\tau=10^{-3}$ after the global refinement $H=\frac14\to h=\frac18$ in \cref{ex:1d-transfers}, with interpolation (left), the ERT (middle) and the ERT-f (right) as mesh transfer.
  The top row shows $u_h^{N+1}$ (solid) and $u(t_{N+1})$ (dashed), and the bottom row shows $d_h^{N+1}$ (solid) and $\partial_tu(t_{N+1})$ (dashed).
  The vertical grid lines mark the nodes of the old mesh $\T_H$.
  }
  \label{fig:example-1d}
\end{figure}
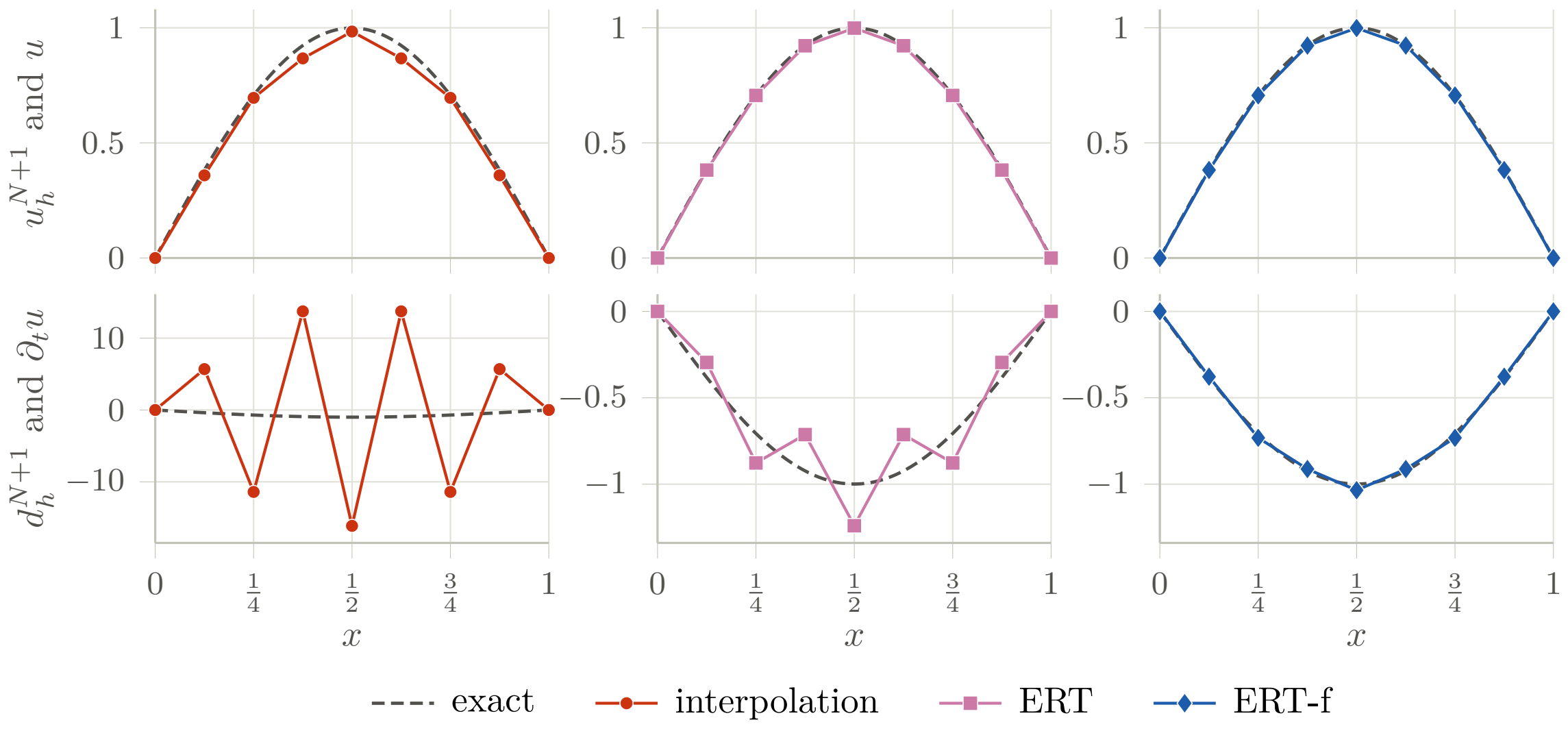

\begin{example}\label{ex:1d-transfers}
    We illustrate the effect of the mesh transfer by comparing interpolation, the ERT and the ERT-f after a global refinement in one spatial dimension on a simple example.
    Let $\Omega = (0, 1)$, $p=1$  and $u(x,t)=e^{t_N-t}\sin(\pi x)$ with $f(x,t)=(\pi^2-1)u(x,t)$, so that $u(t_N)=\sin(\pi x)$.
    On the uniform mesh $\T_H$ with $H=\frac14$,~i.e.,~with five nodes, let $u_H^N := I_H u(t_N)$ be the discrete solution at $t_N$ and $d_H^N:=P_Hf^N-A_Hu_H^N$ the corresponding discrete time derivative according to \eqref{eq:euler}.
    We refine $\T_H$ globally to $\T_h$ with $h=\frac18$, transfer $u_H^N$ by $\Pi\in\{I_h,\PiER,\PiERf\}$ and compute one step of \eqref{eq:intro-transfer} with $\tau_{N+1}=\tau_N=\tau=10^{-3}$.

    \cref{fig:example-1d} shows that, for all three transfers, $u_h^{N+1}$ differs only slightly from $u(t_{N+1})$, while $d_h^{N+1} = (u_h^{N+1} - \Pi u_H^N) / \tau $ depends strongly on the transfer.
    For interpolation, $d_h^{N+1}$ oscillates strongly between the old and the new nodes.
    For the ERT, the residual does not contain the defect $\delta_N$,~cf.~\eqref{eq:comp-interpol-ert}, and thus these oscillations are damped.
    For the ERT-f, these oscillations in $d_h^{N+1}$ are almost completely avoided, since $r_N=0$ by \cref{thm:ERf}.
\end{example}

\subsection{\texorpdfstring{Constrained elliptic reconstruction transfer}{Constrained elliptic reconstruction transfer}}\label{sec:constrained-ell-rec-trans}

If a transfer is required to preserve the old nodal values, only the new degrees of freedom can be determined. 
For instance, nodal interpolation has this property.
We now consider a constrained elliptic reconstruction transfer that also preserves the old nodal values.
Let $V_H\subset V_h$ and define
\begin{align*}
  V_h^\bullet:=\{v\in V_h:\ v(z)=0\ \ \forall z\in\N_H\},
\end{align*}
the hierarchical complement of $V_H$, so that $V_h=V_H\oplus V_h^\bullet$.
For the discrete time derivative transfer, we again choose $\Pi_d=P_h$.

\begin{definition}\label{def:constrained}%
  We define the constrained elliptic reconstruction transfer with data $f$ (CERT-f) by $\PiC u_H^N:=u_H^N+u_h^\bullet$, where $u_h^\bullet\in V_h^\bullet$ satisfies $a(u_H^N+u_h^\bullet,v)=(f^N-d_H^N,v)$ for all $v\in V_h^\bullet$.
\end{definition}

This construction is the Galerkin approximation of the ERT-f problem in the affine space $u_H^N + V_h^\bullet$.
The coarse nodal values serve as Dirichlet data, and we compute only the new degrees of freedom. 
Consequently, we have $I_H\PiC u_H^N=u_H^N$ under nested refinement.

\begin{lemma}\label{lem:constrained}
  Let $\hat u_h^N:=\PiERf u_H^N$ and $z:=\hat u_h^N-\PiC u_H^N$. Then, the following holds.
  \begin{enumerate}[label=(\roman*)]
    \item $a(z,v)=0$ for all $v\in V_h^\bullet$, and $I_Hz=I_H\hat u_h^N-u_H^N$.
    \item The residual is $r_N=A_hz$, and $r_N\perp V_h^\bullet$.
    \item It holds $r_N=0$ if and only if $z=0$, and $z=0$ if and only if $I_H\hat u_h^N=u_H^N$.
  \end{enumerate}
\end{lemma}
\begin{proof}
    By the definition of $\PiERf$ in \cref{def:ERf} and $\PiC$ in \cref{def:constrained}, we have $a(z,v)=(f^N-d_H^N,v)-(f^N-d_H^N,v)=0$ for all $v\in V_h^\bullet$.
    Since the coarse nodal values are preserved, we have $I_Hz=I_H\hat u_h^N-u_H^N$.
    Using $A_h\hat u_h^N=P_h(f^N-d_H^N)$, the transfer residual satisfies
    \begin{align*}
        r_N = P_hf^N-P_h d_H^N-A_h\PiC u_H^N
        = P_h(f^N-d_H^N)-A_h\PiC u_H^N
        =A_h\hat u_h^N-A_h\PiC u_H^N=A_hz.
    \end{align*}
    Finally, we have $(A_hz,v)=a(z,v)=0$ for $v\in V_h^\bullet$.
    For (iii), we note that $A_h$ is injective on $V_h$, since $(A_hv,v)=\norm{\nabla v}^2$ for $v\in V_h$, and thus $r_N=A_hz=0$ if and only if $z=0$.
    If $z=0$, then $I_H\hat u_h^N-u_H^N=I_Hz=0$ by (i).
    Conversely, if $I_H\hat u_h^N=u_H^N$, then $I_Hz=0$ by (i),~i.e.,~$z\in V_h^\bullet$, and (i) yields $\norm{\nabla z}^2=a(z,z)=0$.
\end{proof}

By \cref{lem:constrained}~(iii), the CERT-f yields $r_N=0$ if and only if the ERT-f preserves the old nodal values.
In one spatial dimension for linear finite elements, the CERT-f coincides with the ERT-f.

\begin{proposition}\label{prop:constrained-1d}
  Let $m=1$, $p=1$ and $V_H\subset V_h$, and let $u_H^N,d_H^N\in V_H$ satisfy \eqref{eq:euler} on $\T_H$ for $n=N$.
  Then $\PiC u_H^N=\PiERf u_H^N$ and $r_N=0$.
\end{proposition}
\begin{proof}
  Let $\hat u_h^N:=\PiERf u_H^N$.
  For $v\in V_H\subset V_h$, \cref{def:ERf} and \eqref{eq:euler} on $\T_H$ yield
  \begin{align*}
    a(\hat u_h^N,v)=(f^N-d_H^N,v)=(f^N-P_Hf^N+A_Hu_H^N,v)=(A_Hu_H^N,v)=a(u_H^N,v),
  \end{align*}
  and thus $R_H\hat u_h^N=u_H^N$.
  Since $m=1$ and $p=1$, for $w\in H^1_0(\Omega)$ and $v\in V_H$, the derivative $v'$ is constant on each $K\in\T_H$, and $w-I_Hw$ vanishes at the endpoints of $K$, so that
  \begin{align*}
    a(w-I_Hw,v)=\sum_{K\in\T_H}v'|_K\int_K(w-I_Hw)'\,dx=0.
  \end{align*}
  Thus, $I_Hw=R_Hw$ for all $w\in H^1_0(\Omega)$ by the definition of the Ritz projection.
  Hence, $I_H\hat u_h^N=R_H\hat u_h^N=u_H^N$, and \cref{lem:constrained}~(iii) yields $z=0$,~i.e.,~$\PiC u_H^N=\PiERf u_H^N$, and $r_N=0$.
\end{proof}

\subsection{Summary}

\begin{table}[htb!]
  \caption{
    Comparison of the mesh transfers for nested refinement $V_H\subset V_h$.
    For $\PiC$, the residual is given in terms of $z=\PiERf u_H^N-\PiC u_H^N$,~cf.~\cref{lem:constrained}.
    The column nodal indicates whether the old nodal values are preserved in general,~i.e.,~whether $I_H\Pi u_H^N=u_H^N$.
    The cost Poisson means one solve with the stiffness matrix of $V_h$.
    }
  \label{tab:summary}
  \centering
  \begin{tabular}{@{}l l c l@{}}    
    \toprule
    transfer & residual $r_N$ & nodal & cost \\
    \midrule
    $I_h$ & $(P_h-P_H)f^N-(\Id-P_H)A_hu_H^N$ & yes & (almost) none \\
    $\PiER$ & $(P_h-P_H)f^N$ & no & Poisson \\
    $\PiERf$ & $0$ & no & Poisson \\
    $\PiC$ & $A_hz$ & yes & Poisson (new DoFs) \\
    \bottomrule
  \end{tabular}

\end{table}

\Cref{tab:summary} compares the transfers for nested refinement.
Interpolation, which coincides with the $L^2$ and the Ritz projection on $V_H$, preserves $u_H^N$, but its residual contains the defect $\delta_N$, which in general does not converge for $p=1$,~cf.~\cref{lem:projection}.
The ERT removes $\delta_N$ from the residual and leaves the data oscillation $(P_h-P_H)f^N$,~cf.~\cref{prop:ER}.
The ERT-f yields $r_N=0$ even for arbitrary mesh changes,~cf.~\cref{thm:ERf}, but, like the ERT, it changes the old nodal values in general.
The CERT-f preserves the old nodal values, and its residual vanishes if and only if the ERT-f preserves them as well,~cf.~\cref{lem:constrained}, e.g., for $m=1$ and $p=1$,~cf.~\cref{prop:constrained-1d}.

\section{Space and time adaptive algorithm}\label{sec:algorithm}

In this section, we describe the space and time adaptive algorithm used later in~\cref{sec:num-adaptive}.
It combines the step size control based on the estimator~\eqref{eq:est-euler} with spatial mesh adaptation.

We select the step sizes by the step size control of~\cite[Section~II.4]{HNW93}, a modified version of~\eqref{eq:controller-intro}, based on the estimator~\eqref{eq:est-euler} and a prescribed tolerance $\tol>0$.
Given $u_X^{n-1}$ and $d_X^{n-1}$ on a mesh $\T_X$ at $t_{n-1}$ and a proposed step size $\tau_n$, we compute $u_X^n$ and $d_X^n$ by~\eqref{eq:euler}, the estimator $\eta^n$ by~\eqref{eq:est-euler}, and
\begin{equation}\label{eq:controller}
  \tau_{\mathrm{new}}:=\tau_n\min\bigl\{q_{\max},\max\{0.1,\,0.9\,(\tol/\eta^n)^{1/2}\}\bigr\}.
\end{equation}
If $\eta^n\le\tol$, the step is accepted,~i.e.,~$t_n:=t_{n-1}+\tau_n$, and $\tau_{n+1}:=\tau_{\mathrm{new}}$ is the proposed step size for the next step.
If $\eta^n>\tol$, the step is rejected,~i.e.,~$u_X^n$ is discarded, and step $n$ is recomputed from $u_X^{n-1}$ and $d_X^{n-1}$ with $\tau_n$ replaced by $\tau_{\mathrm{new}}$, until it is accepted.
We set $q_{\max}=2$, and $q_{\max}=1$ if step $n$ has been rejected at least once.
The last step is shortened such that it ends at $T$.

As spatial error indicator, we use the jumps of the normal derivative of the discrete solution across the cell boundaries.
For a cell $K$ of the current mesh $\T_X$ with diameter $h_K$ and the discrete solution $u_X^n$, we define
\begin{equation}\label{eq:indicator}
  \eta_K^2:=\frac{h_K}{24}\sum_{F\subset\partial K}\norm{\jump{\partial_{n_F}u_X^n}_F}_{L^2(F)}^2,
\end{equation}
where the sum runs over all faces $F$ in the interior of $\Omega$, $n_F$ is a unit normal vector of $F$, and $\jump{\partial_{n_F}u_X^n}_F$ denotes the jump of the normal derivative of $u_X^n$ across $F$.
This is the jump term of the residual \emph{a~posteriori} error estimator for the Poisson problem~\cite{KGZB83}, as computed by the Kelly error estimator of deal.II~\cite{ABB26}.
We mark the cells by the bulk criterion~\cite{Doe96}, where the cells with the largest indicators are marked for refinement until their indicators sum up to $60\,\%$ of $\sum_K\eta_K$, and the cells with the smallest indicators for coarsening until their indicators sum up to $40\,\%$ of $\sum_K\eta_K$.
Cells of a prescribed maximal level are not refined, and cells of a prescribed minimal level are not coarsened.

After every fifth accepted time step, we refine and coarsen the marked cells of the current mesh $\T_H:=\T_X$, which yields the new mesh $\T_h$, and continue with $\T_X:=\T_h$.
For all transfers, the estimator and the data of the reconstructions use the $L^2$ projection $P_hd_H^N$ of the discrete time derivative $d_H^N$,~i.e.,~$\Pi_d=P_h$.
If the mesh is changed after step $N$, then $u_h^{N+1}$ is computed by~\eqref{eq:intro-transfer} and $\eta^{N+1}$ by~\eqref{EQ:est-transfer}.

As initial value, we use the Ritz projection $u_X^0:=R_Xu_0$.
The initial mesh $\T_X$ is obtained from a uniform mesh by several adaptations based on the indicator~\eqref{eq:indicator} for $u_X^0$, which is recomputed on every intermediate mesh.
Further, we set $d_X^0=P_X(f^0+\Delta u_0)$.
Then, we choose the first step size $\tau_1:=0.8\,\tol/\norm{d_X^0}$.
\Cref{alg:adaptive} summarizes the space and time adaptive method.

\begin{algorithm}[htb!]
  \caption{Space and time adaptive implicit Euler method with the mesh transfer $\Pi$.}
  \label{alg:adaptive}
  \begin{algorithmic}[1]
    \State Adapt the initial mesh $\T_X$ to $u_0$, compute $u_X^0:=R_Xu_0$ and $d_X^0:=P_Xf^0-A_Xu_X^0$, and set $t_0:=0$, $\tau_1:=0.8\,\tol/\norm{d_X^0}$, and $n:=1$.
    \While{$t_{n-1}<T$}
      \Repeat
        \State Compute $u_X^n$ and $d_X^n$ on the current mesh $\T_X$ by~\eqref{eq:euler}, or by~\eqref{eq:intro-transfer} in the first step after a mesh change.
        \State Compute $\eta^n$ by~\eqref{eq:est-euler}, or by~\eqref{EQ:est-transfer} after a mesh change, and $\tau_{\mathrm{new}}$ by~\eqref{eq:controller}.
        \If{$\eta^n>\tol$}
          \State $\tau_n\gets\tau_{\mathrm{new}}$ \Comment{rejection}
        \EndIf
      \Until{$\eta^n\le\tol$}
      \State $t_n:=t_{n-1}+\tau_n$ and $\tau_{n+1}:=\tau_{\mathrm{new}}$.
      \If{$n$ is a multiple of $5$ and $t_n<T$}
        \State Compute the indicators~\eqref{eq:indicator} for $u_X^n$, mark the cells, and refine and coarsen the mesh $\T_H:=\T_X$ to a new mesh $\T_h$.
        \State Transfer $u_H^N$ by $\Pi$ and $d_H^N$ by $P_h$ to $V_h$, and set $\T_X:=\T_h$.
      \EndIf
      \State $n\gets n+1$
    \EndWhile
  \end{algorithmic}
\end{algorithm}

\section{Numerical experiments}\label{sec:numerics}

In \cref{sec:num-heat}, we use constant step sizes and a single prescribed refinement to isolate the effect of the mesh transfer, and in \cref{sec:num-adaptive}, we apply the space and time adaptive algorithm introduced in \cref{sec:algorithm}.
All computations use the finite element library deal.II~9.8.0~\cite{ABB26} with Lagrange elements $Q_p$ on quadrilaterals.

\subsection{One global refinement}\label{sec:num-heat}

We continue the example of \cref{sec:issue}.
In each run, we compute $40$ steps with a constant step size $\tau\in\{10^{-1},3\cdot10^{-2},10^{-2},\dots,3\cdot10^{-5},10^{-5}\}$.
We compute the estimator $\eta^n$ by~\eqref{eq:est-euler}, and by~\eqref{EQ:est-transfer} for $n=N+1$, and compare it with the leading term $\eta^n_{\mathrm{ex}}:=\frac{\tau^2}{2}\norm{\partial_{tt}u(t_n)}$ of the LTE~\eqref{eq:euler-ex-est}.
For the discrete time derivative, we use $\Pi_d=P_h$, and for $\Pi$, we compare interpolation, the ERT, the ERT-f and the CERT-f,~cf.~\cref{sec:transfers}.
Under nested refinement, the $L^2$ projection and the Ritz projection coincide with interpolation,~cf.~\cref{lem:nesting}.

\begin{figure}[htb!]
  \centering
  \begin{tikzpicture}
  \definecolor{figink}{HTML}{52514e}
  \definecolor{figmuted}{HTML}{898781}
  \definecolor{figgrid}{HTML}{e1e0d9}
  \definecolor{figaxis}{HTML}{c3c2b7}
  \definecolor{figinterp}{HTML}{cc3311}
  \definecolor{figltwo}{HTML}{c98500}
  \definecolor{figert}{HTML}{cc79a7}
  \definecolor{figertf}{HTML}{1c5cab}
  \definecolor{figcertf}{HTML}{648fff}
  \pgfplotsset{
    sfaxis/.style={
      scale only axis,
      axis lines*=left, axis line style={-, figaxis, line width=0.6pt},
      tick align=outside, tick style={figaxis},
      minor tick length=0pt, %
      grid=major, major grid style={figgrid, line width=0.4pt},
      tick label style={font=\scriptsize, text=figink},
      label style={font=\footnotesize, text=figink},
      title style={font=\footnotesize, text=figink, yshift=-0.5ex},
      legend cell align=left,
    },
    sfline/.style={line width=0.9pt, line cap=round, line join=round},
    sfinterp/.style={sfline, figinterp, mark=*, mark size=1.8pt,
      mark options={solid, fill=figinterp, draw=white, line width=0.3pt}},
    sfltwo/.style={sfline, figltwo, dash pattern=on 3pt off 1.5pt, mark=*, mark size=1.3pt,
      mark options={solid, fill=figltwo, draw=white, line width=0.3pt}},
    sfert/.style={sfline, figert, mark=square*, mark size=2pt,
      mark options={solid, fill=figert, draw=white, line width=0.3pt}},
    sfertf/.style={sfline, figertf, mark=diamond*, mark size=2.6pt,
      mark options={solid, fill=figertf, draw=white, line width=0.3pt}},
    sfconstr/.style={sfline, figcertf, mark=triangle*, mark size=2.8pt,
      mark options={solid, fill=figcertf, draw=white, line width=0.3pt, rotate=180}},
  }
  \tikzset{
    sfslope/.style={figmuted, line width=0.5pt},
    sfnote/.style={font=\scriptsize, text=figink},
    sfreference/.style={figink, line width=0.8pt, dash pattern=on 3.5pt off 1.75pt},
  }
  \begin{axis}[sfaxis, name=a, width=3.6cm, height=3.6cm, ymode=log,
      xmin=1, xmax=30, xtick={1,5,10,15,20,25,30}, ymin=3e-2, ymax=1e6,
      ytick={1e-1,1e0,1e1,1e2,1e3,1e4,1e5,1e6},
      xlabel={step $n$}, ylabel={$\eta^n/\eta^n_{\mathrm{ex}}$},
      legend to name=heatlegend, legend columns=4,
      legend style={draw=none, fill=none, font=\footnotesize,
        /tikz/every odd column/.append style={column sep=0.2em},
        /tikz/every even column/.append style={column sep=1.2em}}]
    \addlegendimage{sfinterp}\addlegendentry{interpolation}
    \addlegendimage{sfert}\addlegendentry{ERT}
    \addlegendimage{sfertf}\addlegendentry{ERT-f}
    \addlegendimage{sfconstr}\addlegendentry{CERT-f}
    \draw[figmuted, line width=0.7pt]
      ({axis cs:10.5,1}|-{rel axis cs:0,0}) -- ({axis cs:10.5,1}|-{rel axis cs:0,1})
      node[anchor=north west, xshift=2pt, yshift=-2pt, inner sep=0pt, sfnote] {refined};
    \addplot[sfinterp, no markers] table[row sep=\\, x=x, y=y] {
      x y\\
      2 1.06231\\
      3 1.06219\\
      4 1.06207\\
      5 1.06194\\
      6 1.06182\\
      7 1.0617\\
      8 1.06158\\
      9 1.06146\\
      10 1.06133\\
      11 108322\\
      12 59565.3\\
      13 26822.3\\
      14 12078.1\\
      15 5438.81\\
      16 2449.11\\
      17 1102.84\\
      18 496.61\\
      19 223.624\\
      20 100.698\\
      21 45.3448\\
      22 20.419\\
      23 9.19512\\
      24 4.14151\\
      25 1.86707\\
      26 0.845649\\
      27 0.391853\\
      28 0.199891\\
      29 0.131172\\
      30 0.113489\\
    };
    \addplot[sfert, no markers] table[row sep=\\, x=x, y=y] {
      x y\\
      2 1.06231\\
      3 1.06219\\
      4 1.06207\\
      5 1.06194\\
      6 1.06182\\
      7 1.0617\\
      8 1.06158\\
      9 1.06146\\
      10 1.06133\\
      11 155.957\\
      12 96.6144\\
      13 43.515\\
      14 19.6157\\
      15 8.87909\\
      16 4.09911\\
      17 2.05514\\
      18 1.29339\\
      19 1.07499\\
      20 1.02503\\
      21 1.01458\\
      22 1.01243\\
      23 1.01198\\
      24 1.01186\\
      25 1.01182\\
      26 1.0118\\
      27 1.01177\\
      28 1.01175\\
      29 1.01172\\
      30 1.0117\\
    };
    \addplot[sfertf, no markers] table[row sep=\\, x=x, y=y] {
      x y\\
      2 1.06231\\
      3 1.06219\\
      4 1.06207\\
      5 1.06194\\
      6 1.06182\\
      7 1.0617\\
      8 1.06158\\
      9 1.06146\\
      10 1.06133\\
      11 10.9154\\
      12 4.99784\\
      13 2.42567\\
      14 1.41845\\
      15 1.10763\\
      16 1.03325\\
      17 1.01748\\
      18 1.01423\\
      19 1.01355\\
      20 1.01339\\
      21 1.01334\\
      22 1.0133\\
      23 1.01328\\
      24 1.01325\\
      25 1.01322\\
      26 1.0132\\
      27 1.01317\\
      28 1.01315\\
      29 1.01312\\
      30 1.01309\\
    };
    \addplot[sfconstr, no markers] table[row sep=\\, x=x, y=y] {
      x y\\
      2 1.06231\\
      3 1.06219\\
      4 1.06207\\
      5 1.06194\\
      6 1.06182\\
      7 1.0617\\
      8 1.06158\\
      9 1.06146\\
      10 1.06133\\
      11 797.295\\
      12 407.022\\
      13 154.285\\
      14 63.2622\\
      15 27.286\\
      16 12.1662\\
      17 5.69825\\
      18 3.08102\\
      19 2.2061\\
      20 1.98163\\
      21 1.93165\\
      22 1.91999\\
      23 1.91625\\
      24 1.91412\\
      25 1.91233\\
      26 1.9106\\
      27 1.90889\\
      28 1.90719\\
      29 1.90549\\
      30 1.90379\\
    };
  \end{axis}
  \begin{axis}[sfaxis, name=b, at={(a.south east)}, anchor=south west, xshift=2.07cm,
      width=3.6cm, height=3.6cm, xmode=log, ymode=log, x dir=reverse,
      xmin=5e-6, xmax=0.2, xtick={1e-5,1e-4,1e-3,1e-2,1e-1}, ymin=0.03, ymax=1e+07, ytick={1e-1,1e0,1e1,1e2,1e3,1e4,1e5,1e6,1e7},
      xlabel={time step $\tau$}, ylabel={$\eta^{N+1}/\eta^{N+1}_{\mathrm{ex}}$}]
    \draw[sfslope] (axis cs:0.1,100000) -- (axis cs:0.0316228,316228) -- (axis cs:0.0316228,100000) -- cycle;
    \node[sfnote, anchor=west] at (axis cs:0.0316228,177828) {$\tau^{-1}$};
    \draw[sfslope] (axis cs:0.1,1000) -- (axis cs:0.0316228,10000) -- (axis cs:0.0316228,1000) -- cycle;
    \node[sfnote, anchor=west] at (axis cs:0.0316228,3162.28) {$\tau^{-2}$};
    \addplot[sfinterp] table[row sep=\\, x=x, y=y] {
      x y\\
      0.1 1.28059\\
      0.03 2.97524\\
      0.01 20.1662\\
      0.003 213.702\\
      0.001 1821.8\\
      0.0003 17198.6\\
      0.0001 108322\\
      3e-05 586935\\
      1e-05 2.14399e+06\\
    };
    \addplot[sfert] table[row sep=\\, x=x, y=y] {
      x y\\
      0.1 1.08903\\
      0.03 0.999202\\
      0.01 0.98193\\
      0.003 1.04552\\
      0.001 1.44945\\
      0.0003 21.3502\\
      0.0001 155.957\\
      3e-05 886.164\\
      1e-05 3279.98\\
    };
    \addplot[sfertf] table[row sep=\\, x=x, y=y] {
      x y\\
      0.1 1.08876\\
      0.03 0.997858\\
      0.01 0.981534\\
      0.003 1.16723\\
      0.001 2.08232\\
      0.0003 5.27455\\
      0.0001 10.9154\\
      3e-05 17.6954\\
      1e-05 21.5347\\
    };
    \addplot[sfconstr] table[row sep=\\, x=x, y=y] {
      x y\\
      0.1 0.914476\\
      0.03 0.0600369\\
      0.01 2.98939\\
      0.003 13.8187\\
      0.001 45.7972\\
      0.0003 180.961\\
      0.0001 797.295\\
      3e-05 4161.38\\
      1e-05 15843.1\\
    };
  \end{axis}
  \begin{axis}[sfaxis, name=c, at={(b.south east)}, anchor=south west, xshift=1.18cm,
      width=3.6cm, height=3.6cm, xmode=log, ymode=log, x dir=reverse,
      xmin=5e-6, xmax=0.2, xtick={1e-5,1e-4,1e-3,1e-2,1e-1}, ymin=0.3, ymax=300000, ytick={1e0,1e1,1e2,1e3,1e4,1e5},
      xlabel={time step $\tau$}]
    \draw[sfslope] (axis cs:0.1,300) -- (axis cs:0.01,3000) -- (axis cs:0.01,300) -- cycle;
    \node[sfnote, anchor=west] at (axis cs:0.01,948.683) {$\tau^{-1}$};
    \addplot[sfinterp] table[row sep=\\, x=x, y=y] {
      x y\\
      0.1 1.10928\\
      0.03 1.03174\\
      0.01 1.16604\\
      0.003 6.63632\\
      0.001 55.5463\\
      0.0003 510.937\\
      0.0001 3085.9\\
      3e-05 15965.9\\
      1e-05 56953.2\\
    };
    \addplot[sfert] table[row sep=\\, x=x, y=y] {
      x y\\
      0.1 1.10912\\
      0.03 1.02871\\
      0.01 0.994265\\
      0.003 0.985986\\
      0.001 0.993625\\
      0.0003 1.26068\\
      0.0001 5.36758\\
      3e-05 28.4222\\
      1e-05 102.4\\
    };
    \addplot[sfertf] table[row sep=\\, x=x, y=y] {
      x y\\
      0.1 1.10912\\
      0.03 1.02871\\
      0.01 0.994265\\
      0.003 0.986131\\
      0.001 0.994098\\
      0.0003 1.00909\\
      0.0001 1.04568\\
      3e-05 1.10859\\
      1e-05 1.15169\\
    };
    \addplot[sfconstr] table[row sep=\\, x=x, y=y] {
      x y\\
      0.1 1.10903\\
      0.03 1.02822\\
      0.01 0.992349\\
      0.003 0.979191\\
      0.001 0.987801\\
      0.0003 1.92189\\
      0.0001 11.2098\\
      3e-05 71.9956\\
      1e-05 334.956\\
    };
  \end{axis}
  \coordinate (labelrow) at (current bounding box.south);
  \node[anchor=north, yshift=-1pt, font=\footnotesize, text=figink] at (a.south |- labelrow) {(a) $p=1$, $\tau=10^{-4}$};
  \node[anchor=north, yshift=-1pt, font=\footnotesize, text=figink] at (b.south |- labelrow) {(b) $p=1$};
  \node[anchor=north, yshift=-1pt, font=\footnotesize, text=figink] at (c.south |- labelrow) {(c) $p=2$};
  \node[anchor=north, yshift=-1mm] at (current bounding box.south) {\pgfplotslegendfromname{heatlegend}};
\end{tikzpicture}
  \vspace{-1cm}
  \caption{
  Heat equation with one global refinement $H=\frac1{16}\to h=\frac1{32}$ after step $N=10$.
  (a) Ratio $\rho^n=\eta^n/\eta^n_{\mathrm{ex}}$ of the estimator and the leading term of the LTE for $p=1$ and $\tau=10^{-4}$.
  (b), (c) Ratio $\rho^{N+1}=\norm{D_h^{N+1}}/\norm{\partial_{tt}u(t_{N+1})}$ in the first step after the refinement against the constant step size $\tau$ for $p=1$ and $p=2$.
  }
  \label{fig:heat}
\end{figure}
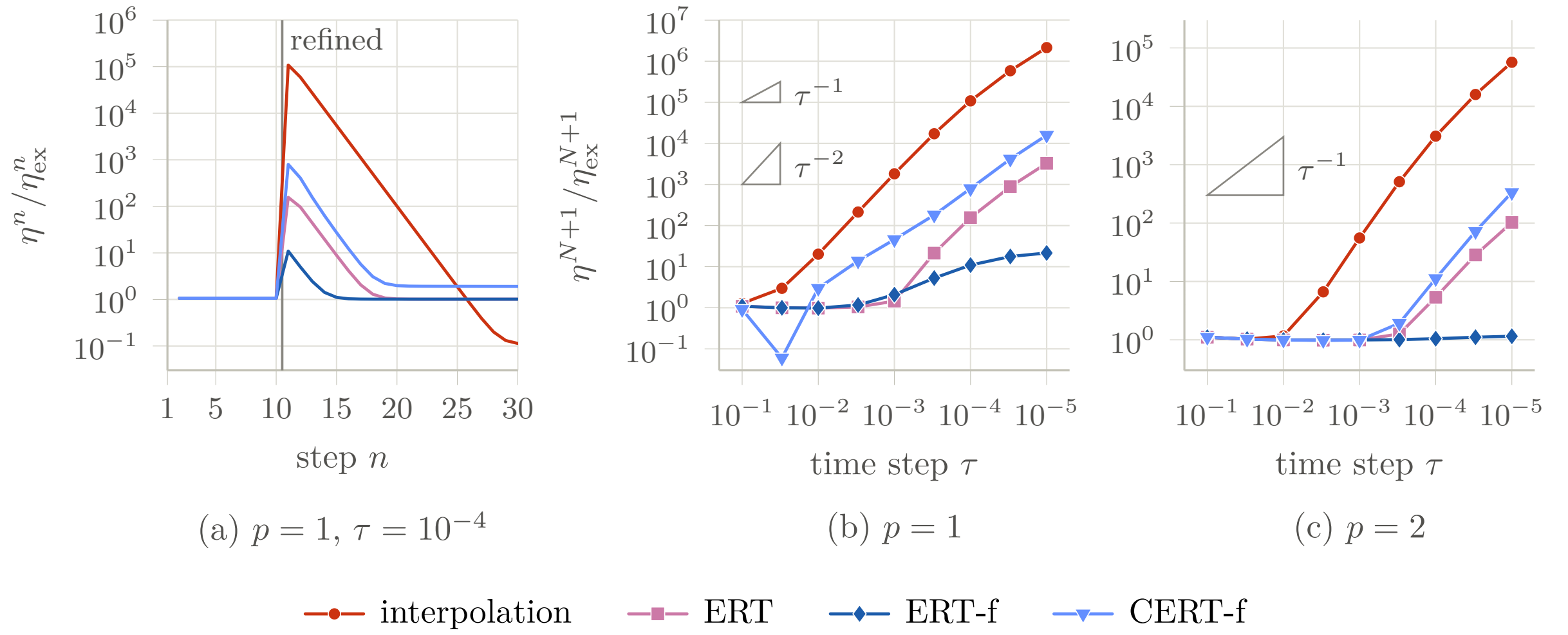

\Cref{fig:heat}(a) shows the ratio $\rho^n:=\eta^n/\eta^n_{\mathrm{ex}}$ for $p=1$ and $\tau=10^{-4}$.
Before the refinement, the estimator agrees with the local error,~i.e.,~$\rho^n\approx1$.
In the first step after the refinement, the ratio increases by about five orders of magnitude for interpolation, by two to three orders for the CERT-f and the ERT, and by one order for the ERT-f.
For the ERT and the ERT-f, it returns to one within a few steps.
For interpolation and the CERT-f, it does not return to one within the computed steps, and for interpolation, the estimator even underestimates the local error by about one order of magnitude.

\Cref{fig:heat}(b) and~(c) show the ratio $\rho^{N+1}$ in the first step after the refinement against $\tau$ for $p=1$ and $p=2$.
For constant step sizes, $\rho^{N+1}=\norm{D_h^{N+1}}/\norm{\partial_{tt}u(t_{N+1})}$ is the factor by which the estimator overestimates the local error.
Unless $r_N=0$, this factor grows like $\tau^{-1}$ for small $\tau$,~cf.~\cref{thm:overestimation}.
For larger $\tau$, it grows like $\tau^{-2}$ for interpolation, since the contribution of $r_N$ is damped additionally according to \cref{lem:damping}.
For the ERT, the residual is only the data oscillation,~cf.~\cref{prop:ER}, and the overestimation sets in only for smaller step sizes.
For the ERT-f, the ratio remains bounded, since $r_N=0$ by \cref{thm:ERf}, and its moderate increase is caused by the defect $\delta_d$ of the transferred time derivative in~\eqref{eq:dn1}.
For $p=2$, the overestimation is smaller for all transfers.

\subsection{Moving bump with adaptivity in space and time}\label{sec:num-adaptive}

\begin{figure}[htb!]
  \centering
  \begin{minipage}[t]{0.3\textwidth}
    \centering
    \includegraphics[width=\textwidth]{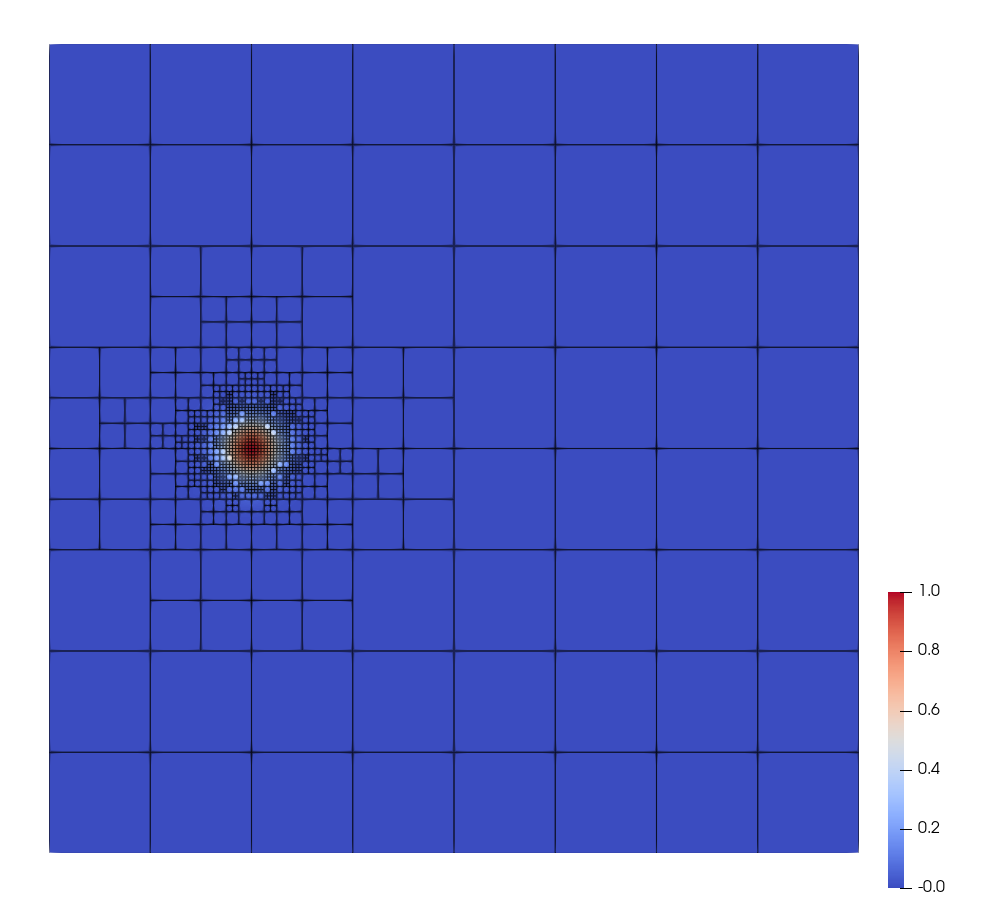}\\
    {\footnotesize (a) $t\approx0.25$}
  \end{minipage}
  \begin{minipage}[t]{0.3\textwidth}
    \centering
    \includegraphics[width=\textwidth]{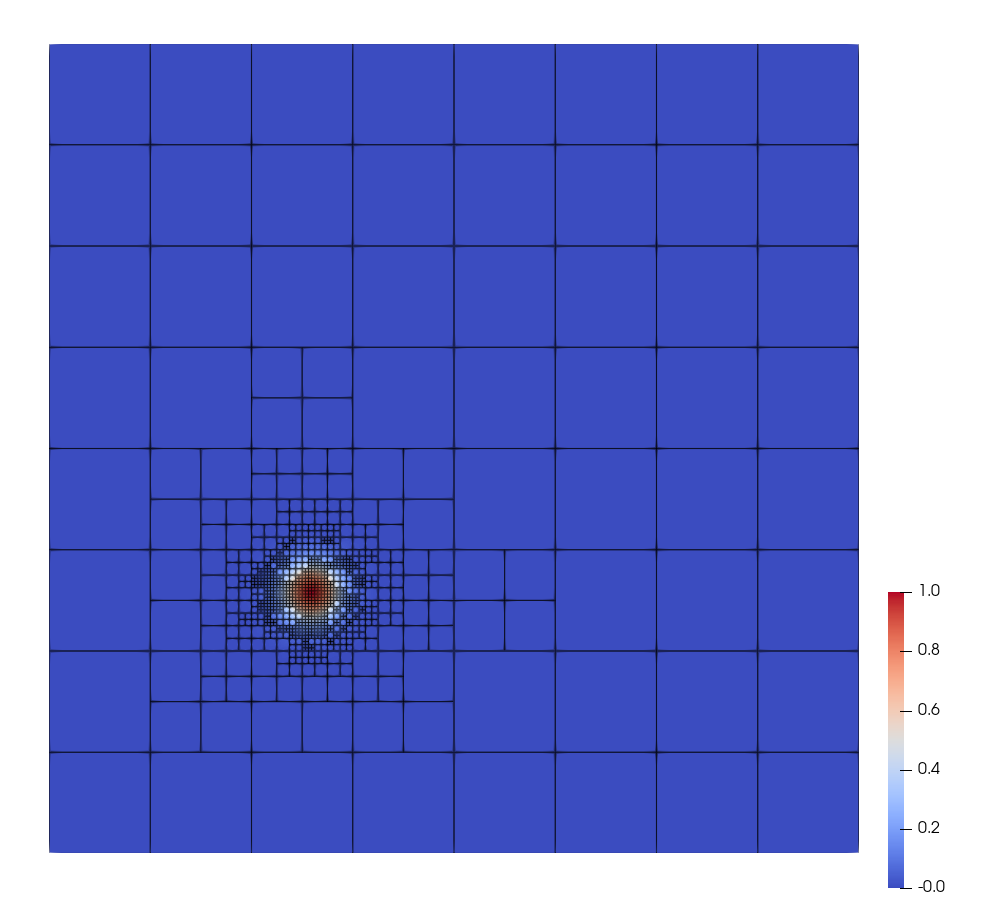}\\
    {\footnotesize (b) $t\approx0.50$}
  \end{minipage}
  \begin{minipage}[t]{0.3\textwidth}
    \centering
    \includegraphics[width=\textwidth]{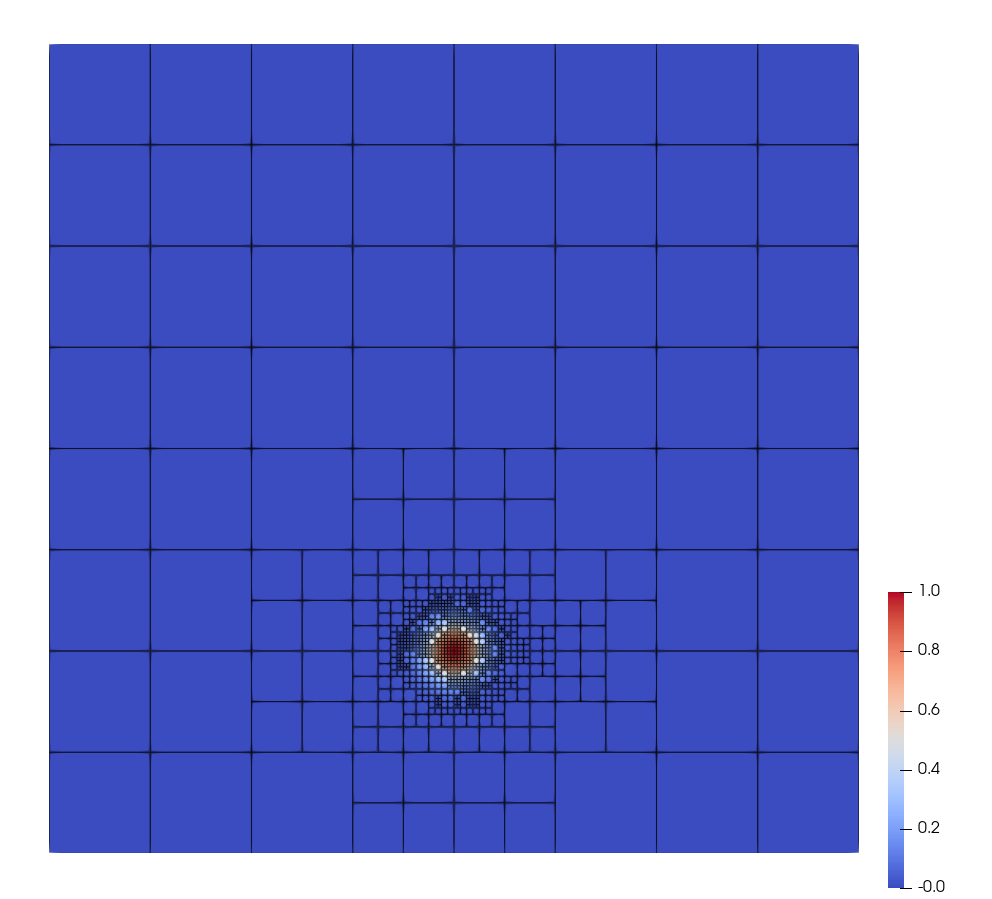}\\
    {\footnotesize (c) $t\approx0.75$}
  \end{minipage}
  \caption{Moving bump: Snapshots of adaptive meshes and the solution using ERT-f with $\tol=10^{-4}$.}
  \label{fig:snapshots}
\end{figure}

Next, we apply the space and time adaptive algorithm introduced in \cref{sec:algorithm} and demonstrate the influence of the mesh transfer.
We consider the domain $\Omega=(-1,1)^2$, final time $T=1$, and
\begin{align*}
  u(x,t)=\exp\Bigl(-\frac{\abs{x-c(t)}^2}{2\sigma^2}\Bigr),\quad
  c(t)=\tfrac12\bigl(\cos\theta(t),\sin\theta(t)\bigr),\quad
  \theta(t)=\tfrac34\pi+\pi t,\quad\sigma=0.05,
\end{align*}
with $f(x, t):=\partial_tu(x, t)-\Delta u(x,t)$.
The bump moves on a half circle of radius $\frac12$.
We impose homogeneous Dirichlet boundary conditions.
We use $p=1$.
Cells of level $8$ are not refined, and cells of level $3$ are not coarsened, so that the cell sizes are between $\frac1{128}$ and $\frac14$.
The initial uniform mesh with mesh size $h=\frac14$ is adapted five times to $u(0)$.
Snapshots of the solution and the adapted mesh are shown in \cref{fig:snapshots}.
We compare interpolation, the $L^2$ projection, the ERT, the ERT-f and the CERT-f.
We apply \cref{alg:adaptive} with the tolerances $\tol\in\{10^{-3},3\cdot10^{-4},10^{-4},3\cdot10^{-5},10^{-5}\}$.
For reference, we use fixed uniform meshes with $h=\frac1{64}$ and $h=\frac1{128}$ without mesh changes, at $\tol=10^{-4}$.

\begin{table}[htb!]
  \centering\small
  \caption{
    Moving bump with adaptivity in space and time. Number of attempted time steps, with the number of rejected steps in brackets, and error at the final time $T=1$ against the tolerance. The fixed uniform meshes are only computed for $\tol=10^{-4}$.
    }
  \label{tab:adaptive}
  \medskip
  \begin{tabular}{@{}lrrrrr@{}}
\toprule
tol & $10^{-3}$ & $3\cdot10^{-4}$ & $10^{-4}$ & $3\cdot10^{-5}$ & $10^{-5}$\\
\midrule
\multicolumn{6}{@{}l}{\emph{attempted time steps (rejected steps)}}\\
interpolation & 220\,(1) & 415\,(5) & 1552\,(547) & 143174\,(58944) & 551633\,(157610)\\
$L^2$ projection & 220\,(1) & 417\,(5) & 1921\,(817) & 213585\,(80113) & 754849\,(215672)\\
ERT & 218\,(0) & 397\,(0) & 687\,(0) & 1268\,(4) & 2718\,(333)\\
ERT-f & 218\,(0) & 397\,(0) & 686\,(0) & 1253\,(0) & 2174\,(0)\\
CERT-f & 220\,(1) & 408\,(3) & 1021\,(199) & 27562\,(13259) & 303398\,(113387)\\
fixed mesh $h=1/64$ & -- & -- & 683\,(0) & -- & --\\
fixed mesh $h=1/128$ & -- & -- & 687\,(0) & -- & --\\
\midrule
\multicolumn{6}{@{}l}{\emph{$\norm{u_X^n-u(T)}$ at $T=1$}}\\
interpolation & $6.5\cdot10^{-4}$ & $6.8\cdot10^{-4}$ & $1.2\cdot10^{-3}$ & $7.7\cdot10^{-2}$ & $2.0\cdot10^{-1}$\\
$L^2$ projection & $6.6\cdot10^{-4}$ & $5.3\cdot10^{-4}$ & $4.9\cdot10^{-4}$ & $5.1\cdot10^{-4}$ & $5.9\cdot10^{-4}$\\
ERT & $6.1\cdot10^{-4}$ & $5.0\cdot10^{-4}$ & $4.8\cdot10^{-4}$ & $4.6\cdot10^{-4}$ & $4.9\cdot10^{-4}$\\
ERT-f & $6.1\cdot10^{-4}$ & $5.1\cdot10^{-4}$ & $5.0\cdot10^{-4}$ & $4.6\cdot10^{-4}$ & $4.8\cdot10^{-4}$\\
CERT-f & $6.4\cdot10^{-4}$ & $5.3\cdot10^{-4}$ & $6.5\cdot10^{-4}$ & $6.9\cdot10^{-3}$ & $5.6\cdot10^{-2}$\\
fixed mesh $h=1/64$ & -- & -- & $9.2\cdot10^{-4}$ & -- & --\\
fixed mesh $h=1/128$ & -- & -- & $2.8\cdot10^{-4}$ & -- & --\\
\bottomrule
\end{tabular}

\end{table}

\begin{figure}[htb!]
  \centering
  \begin{tikzpicture}
  \definecolor{figink}{HTML}{52514e}
  \definecolor{figmuted}{HTML}{898781}
  \definecolor{figgrid}{HTML}{e1e0d9}
  \definecolor{figaxis}{HTML}{c3c2b7}
  \definecolor{figinterp}{HTML}{cc3311}
  \definecolor{figltwo}{HTML}{c98500}
  \definecolor{figert}{HTML}{cc79a7}
  \definecolor{figertf}{HTML}{1c5cab}
  \definecolor{figcertf}{HTML}{648fff}
  \pgfplotsset{
    sfaxis/.style={
      scale only axis,
      axis lines*=left, axis line style={-, figaxis, line width=0.6pt},
      tick align=outside, tick style={figaxis},
      minor tick length=0pt, %
      grid=major, major grid style={figgrid, line width=0.4pt},
      tick label style={font=\scriptsize, text=figink},
      label style={font=\footnotesize, text=figink},
      title style={font=\footnotesize, text=figink, yshift=-0.5ex},
      legend cell align=left,
    },
    sfline/.style={line width=0.9pt, line cap=round, line join=round},
    sfinterp/.style={sfline, figinterp, mark=*, mark size=2.5pt,
      mark options={solid, fill=figinterp, draw=white, line width=0.3pt}},
    sfltwo/.style={sfline, figltwo, dash pattern=on 3pt off 1.5pt, mark=*, mark size=2.5pt,
      mark options={solid, fill=figltwo, draw=white, line width=0.3pt}},
    sfert/.style={sfline, figert, mark=square*, mark size=2.5pt,
      mark options={solid, fill=figert, draw=white, line width=0.3pt}},
    sfertf/.style={sfline, figertf, mark=diamond*, mark size=3.2pt,
      mark options={solid, fill=figertf, draw=white, line width=0.3pt}},
    sfconstr/.style={sfline, figcertf, mark=triangle*, mark size=3.2pt,
      mark options={solid, fill=figcertf, draw=white, line width=0.3pt, rotate=180}},
  }
  \pgfplotsset{%
    sfltwo/.append style={solid},
  }
  \tikzset{
    sfslope/.style={figmuted, line width=0.5pt},
    sfnote/.style={font=\scriptsize, text=figink},
    sfreference/.style={figink, line width=0.8pt, dash pattern=on 3.5pt off 1.75pt},
  }
  \begin{axis}[sfaxis, width=6.4cm, height=4.0cm, xmode=log, ymode=log,
      x dir=reverse, xmin=5e-6, xmax=2e-3, ymin=8e-5, ymax=0.5,
      xtick={1e-5,1e-4,1e-3}, ytick={1e-4,1e-3,1e-2,1e-1},
      xlabel={tolerance $\mathrm{tol}$}, ylabel={$\norm{u_X^n-u(T)}$ at $T=1$},
      legend style={draw=none, fill=none, font=\footnotesize, at={(1.04,0.5)}, anchor=west,
        row sep=1pt}]
    \addlegendimage{sfinterp}\addlegendentry{interpolation}
    \addlegendimage{sfltwo}\addlegendentry{$L^2$ projection}
    \addlegendimage{sfert}\addlegendentry{ERT}
    \addlegendimage{sfertf}\addlegendentry{ERT-f}
    \addlegendimage{sfconstr}\addlegendentry{CERT-f}
    \draw[sfreference] (axis cs:2e-3,0.000919072) -- (axis cs:5e-6,0.000919072)
      node[sfnote, anchor=south west, pos=0, inner sep=1.5pt] {fixed mesh $h=1/64$};
    \draw[sfreference] (axis cs:2e-3,0.00027881) -- (axis cs:5e-6,0.00027881)
      node[sfnote, anchor=north east, pos=1, inner sep=1.5pt] {fixed mesh $h=1/128$};
    \addplot[sfinterp] table[row sep=\\, x=x, y=y] {
      x y\\
      0.001 0.000647512\\
      0.0003 0.00067593\\
      0.0001 0.00115475\\
      3e-05 0.0773896\\
      1e-05 0.203755\\
    };
    \addplot[sfltwo] table[row sep=\\, x=x, y=y] {
      x y\\
      0.001 0.000662822\\
      0.0003 0.000529006\\
      0.0001 0.000485343\\
      3e-05 0.000510526\\
      1e-05 0.000586109\\
    };
    \addplot[sfert] table[row sep=\\, x=x, y=y] {
      x y\\
      0.001 0.000613229\\
      0.0003 0.000504282\\
      0.0001 0.000482206\\
      3e-05 0.000457711\\
      1e-05 0.000491284\\
    };
    \addplot[sfertf] table[row sep=\\, x=x, y=y] {
      x y\\
      0.001 0.000612336\\
      0.0003 0.000511443\\
      0.0001 0.000502094\\
      3e-05 0.000456844\\
      1e-05 0.000477412\\
    };
    \addplot[sfconstr] table[row sep=\\, x=x, y=y] {
      x y\\
      0.001 0.000642829\\
      0.0003 0.000529414\\
      0.0001 0.000649457\\
      3e-05 0.00694361\\
      1e-05 0.0558428\\
    };
  \end{axis}
\end{tikzpicture}
  \caption{
    Moving bump example.
    Error at the final time $T=1$ vs.~the time tolerance. 
    The dashed lines show the errors on the fixed uniform meshes for $\tol=10^{-4}$.
    }
  \label{fig:adaptive-cost}
\end{figure}
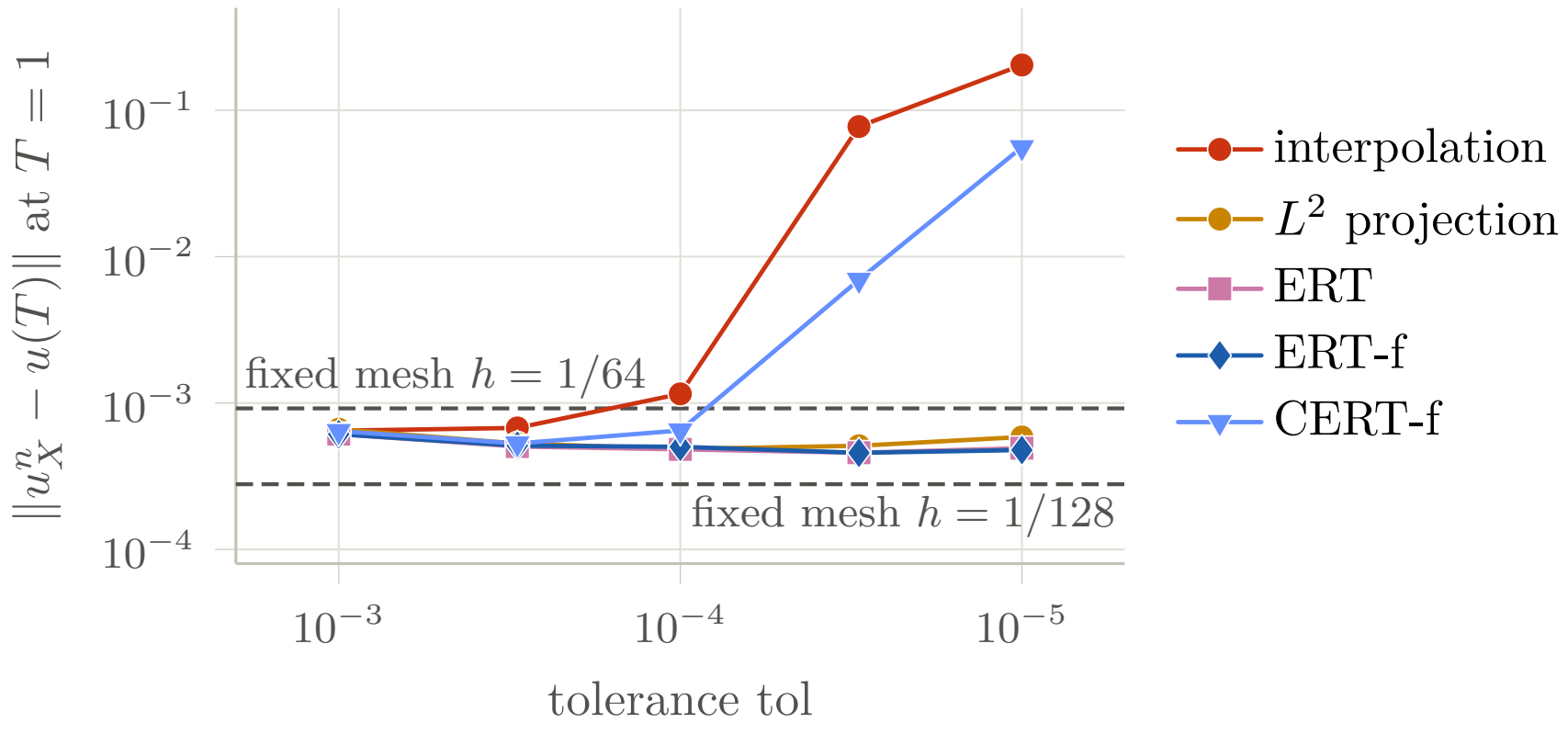

\begin{figure}[htb!]
  \centering
  \input{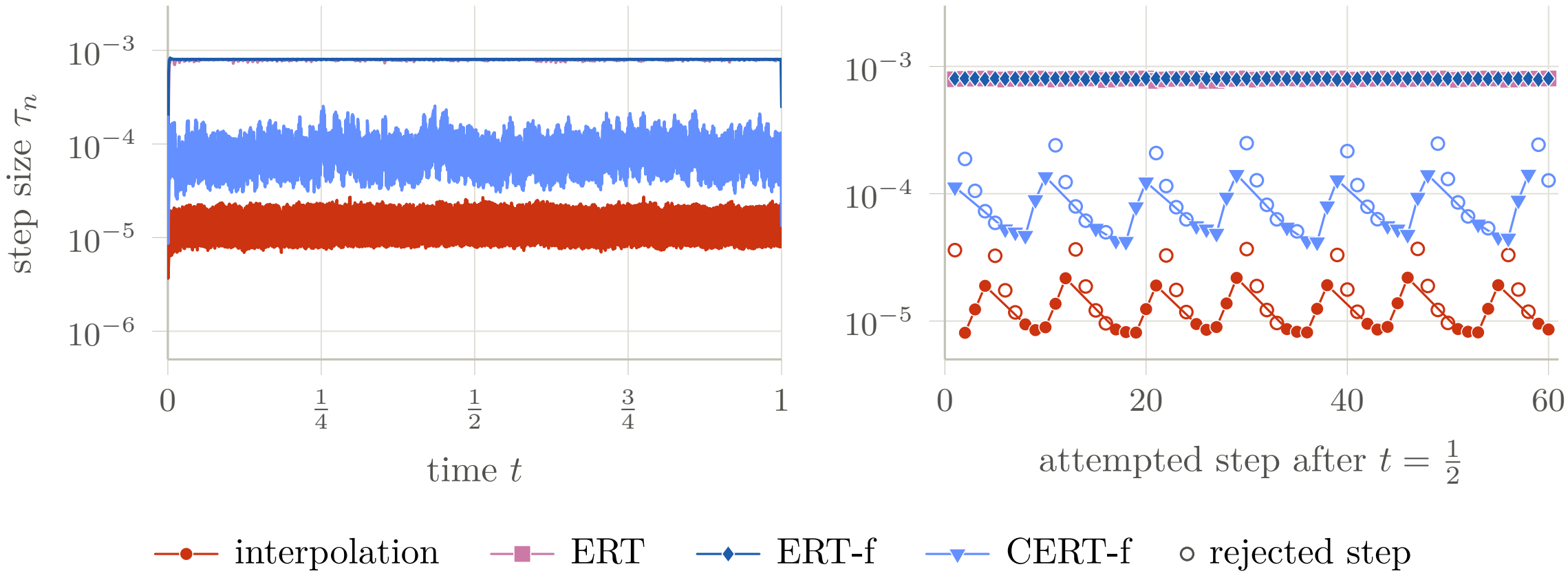}
  \caption{
    Moving bump for $\tol=3\cdot10^{-5}$. 
    Left: Accepted step sizes against time. 
    Right: Step sizes of the first $60$ attempted steps after $t=\frac12$, where rejected steps are marked by circles.
    }
  \label{fig:adaptive-steps}
\end{figure}

\Cref{tab:adaptive} and \cref{fig:adaptive-cost,fig:adaptive-steps} show the results.
For $\tol\ge3\cdot10^{-4}$, all transfers need almost the same number of time steps.
For smaller tolerances, the transfer residual of interpolation and the $L^2$ projection dominates the estimator in the first step after every mesh change,~cf.~\cref{thm:overestimation}.
The step size then oscillates strongly and remains small,~cf.~\cref{fig:adaptive-steps}, so that the number of time steps, and thus of mesh changes, increases much faster than for the ERT-f as the tolerance decreases.
With interpolation, the frequent mesh changes moreover destroy the accuracy, in accordance with~\cite{Dup82}, whereas the $L^2$ projection keeps the error at the level of the spatial discretization.
The CERT-f reduces the transfer residual, but does not remove it,~cf.~\cref{lem:constrained}, and thus behaves similarly to interpolation, with fewer steps and a smaller loss of accuracy.
For the ERT, rejections occur only for the smallest tolerances, caused by the data oscillation in its residual,~cf.~\cref{prop:ER}, and it needs only slightly more steps than the ERT-f.
For the ERT-f, no step is rejected for any tolerance, it needs almost the same number of steps as the fixed uniform meshes, and its error lies between the errors of the two fixed meshes, which have considerably more spatial degrees of freedom.

\section{Conclusions}\label{sec:conclusions}

We have shown that, after a mesh change, the transferred discrete solution and time derivative in general leave a transfer residual in the discrete equation on the new mesh, which enters the discrete time derivative of the first step unchanged,~cf.~\cref{lem:first-step}.
Unless it vanishes, the estimator overestimates the local error by a factor proportional to the inverse step size,~cf.~\cref{thm:overestimation}, in particular after refinement with interpolation.
In adaptive computations, the step size control then reduces the step size after mesh changes, although the local error does not require it.
For small tolerances, this leads to strongly oscillating step sizes.

The ERT-f makes the transfer residual vanish for arbitrary mesh changes at the cost of one elliptic solve,~cf.~\cref{thm:ERf}.
The ERT leaves the data oscillation as residual, so that the overestimation appears only for small step sizes.
The CERT-f, which preserves the old nodal values, reduces the residual, but in general does not remove it,~cf.~\cref{lem:constrained}.

{\sloppy
\subsection*{Acknowledgement}
The author gratefully acknowledges the support of the ERC grant OPTIMAL (grant agreement No. 101125225).
\par}

\bibliographystyle{abbrvurl}
\bibliography{references_new}

\end{document}